\documentclass[9pt,twoside]{article}

\usepackage{amsbsy,amsfonts,amsmath,amssymb,eucal,mathrsfs}
\usepackage[all]{xy}
\usepackage{pstricks,epsfig}
\usepackage{pst-plot}

\usepackage{showlabels}

\newtheorem{definition}{Definition}[section]

\newtheorem{notation}[definition]{Notation}

\newtheorem{construction}[definition]{Construction}

\newtheorem{remark}[definition]{Remark}

\newtheorem{remarks}[definition]{Remarks}

\newtheorem{example}[definition]{Example}

\newtheorem{examples}[definition]{Examples}

\newtheorem{nothing}[definition]{$\!\!$}

\newtheorem{conjecture}[definition]{Conjecture}

\newtheorem{definition*}{Definition}[section]
\newenvironment{defi*}{\begin{definition*} \rm}{\end{definition*}}
\newtheorem{definitions*}[definition*]{Definitions}
\newenvironment{defis*}{\begin{definitions*} \rm}{\end{definitions*}}
\newtheorem{prop*}[definition*]{Proposition}
\newtheorem{lemm*}[definition*]{Lemma}
\newtheorem{coro*}[definition*]{Corollary}
\newtheorem{theo*}[definition*]{Theorem}
\newtheorem{remark*}[definition*]{Remark}
\newenvironment{rema*}{\begin{remark*} \rm}{\end{remark*}}
\newtheorem{remarks*}[definition*]{Remarks}
\newenvironment{remas*}{\begin{remarks*} \rm}{\end{remarks*}}
\newtheorem{example*}[definition*]{Example}
\newenvironment{exam*}{\begin{example*} \rm}{\end{example*}}
\newtheorem{examples*}[definition*]{Examples}
\newenvironment{exams*}{\begin{examples*} \rm}{\end{examples*}}
\newtheorem{nothing*}[definition*]{$\!\!$}
\newenvironment{noth*}{\begin{nothing*} \rm}{\end{nothing*}}

\newtheorem{commentaire*}[definition*]{Commentaire}

\begin{document}

\def \pp {{\odot}}
\def \Rh {{\widehat{R}}}
\def \Sh {{\widehat{S}}}
\def \supp {{\rm Supp}}
\def \codim {{\rm Codim}}
\def \b {{\beta}}
\def \T {{\Theta}}
\def \t {{\theta}}
\def \L {{\cal L}}
\def \sca #1#2{\langle #1,#2 \rangle}
\def\pt{\{{\rm pt}\}}
\def\x {{\underline{x}}}
\def\y {{\underline{y}}}
\def\aut{{\rm Aut}}
\def\ra{\rightarrow}
\def\s{\sigma}\def\OO{\mathbb O}\def\PP{\mathbb P}\def\QQ{\mathbb Q}
 \def\CC{\mathbb C} \def\ZZ{\mathbb Z}\def\JO{{\mathcal J}_3(\OO)}
\newcommand{\G}{\mathbb{G}}
\def\proof{\noindent {\it Proof.}\;}
\def\qed{\hfill $\square$}
\def \uh {{\widehat{u}}}
\def \vh {{\widehat{v}}}
\def \uveeh {{\widehat{u^\vee}}}
\def \vveeh {{\widehat{v^\vee}}}
\def \fh {{\widehat{f}}}
\def \wh {{\widehat{w}}}
\def \wt {{\widetilde{w}}}
\def \Wh {{{W_{{\rm aff}}}}}
\def \Wt {{\widetilde{W}_{{\rm aff}}}}
\def \Qt {{\widetilde{Q}}}
\def \Waff {{W_{{\rm aff}}}}
\def \Waffm {{W_{{\rm aff}}^-}}
\def \Wpaff {{{(W^P)}_{{\rm aff}}}}
\def \Wtpaff {{{(\widetilde{W}^P)}_{{\rm aff}}}}
\def \Wtaffm {{\widetilde{W}_{{\rm aff}}^-}}
\def \lh {{\widehat{\lambda}}}
\def \pit {{\widetilde{\pi}}}
\def \lt {{{\lambda}}}
\def \xh {{\widehat{x}}}
\def \yh {{\widehat{y}}}
\def \a {\alpha}
\def \b {\beta}
\def \l {\lambda}
\def \t {\theta}


\newcommand{\expxy}{\exp_{x \to y}}
\newcommand{\drat}{d_{\rm rat}}
\newcommand{\dmax}{d_{\rm max}}
\newcommand{\zl}{Z(x,L_x,y,L_y)}


\newcommand{\N}{\mathbb{N}}
\newcommand{\A}{{\mathbb{A}_{\rm Aff}}}
\newcommand{\Ah}{{\mathbb{A}_{\rm Aff}}}
\newcommand{\At}{{\widetilde{\mathbb{A}}_{\rm Aff}}}
\newcommand{\Ht}{{{H}^T_*(\Omega K^{\ad})}}
\renewcommand{\H}{{\rm Hi}}
\newcommand{\Ih}{{I_{\rm Aff}}}
\newcommand{\psit}{{\widetilde{\psi}}}
\newcommand{\xit}{{\widetilde{\xi}}}
\newcommand{\Jt}{{\widetilde{J}}}
\newcommand{\Zt}{{\widetilde{Z}}}
\newcommand{\Xt}{{\widetilde{X}}}
\newcommand{\at}{{\widetilde{A}}}
\newcommand{\Z}{\mathbb Z}
\renewcommand{\S}{{\mathbb S}}
\newcommand{\fgcoad}{{G/P_\theta}}
\newcommand{\R}{\mathbb{R}}
\newcommand{\Q}{\mathbb{Q}}
\newcommand{\C}{\mathbb{C}}
\renewcommand{\O}{\mathbb{O}}
\newcommand{\F}{\mathbb{F}}
\newcommand{\p}{\mathbb{P}}
\newcommand{\co}{{\cal O}}
\newcommand{\pos}{{\bf P}}

\renewcommand{\a}{{\alpha}}
\newcommand{\az}{\a_\Z}
\newcommand{\ak}{\a_k}

\newcommand{\rc}{\R_\C}
\newcommand{\cc}{\C_\C}
\newcommand{\hc}{\H_\C}
\newcommand{\oc}{\O_\C}

\newcommand{\rk}{\R_k}
\newcommand{\ck}{\C_k}
\newcommand{\hk}{\H_k}
\newcommand{\ok}{\O_k}

\newcommand{\rz}{\R_Z}
\newcommand{\cz}{\C_Z}
\newcommand{\hz}{\H_Z}
\newcommand{\oz}{\O_Z}

\newcommand{\RR}{\R_R}
\newcommand{\CR}{\C_R}
\newcommand{\HR}{\H_R}
\newcommand{\OR}{\O_R}

\newcommand{\re}{\mathtt{Re}}

\newcommand{\matttr}[9]{
\left (
\begin{array}{ccc}
{} \hspace{-.2cm} #1 & {} \hspace{-.2cm} #2 & {} \hspace{-.2cm} #3 \\
{} \hspace{-.2cm} #4 & {} \hspace{-.2cm} #5 & {} \hspace{-.2cm} #6 \\
{} \hspace{-.2cm} #7 & {} \hspace{-.2cm} #8 & {} \hspace{-.2cm} #9
\end{array}
\hspace{-.15cm}
\right )   }


\newcommand{\dual}{{\bf v}}
\newcommand{\com}{\mathtt{Com}}
\newcommand{\rg}{\mathtt{rg}}
\newcommand{\pu}{{\mathbb{P}^1}}
\newcommand{\scal}[1]{\langle #1 \rangle}
\newcommand{\MK}[2]{{\overline{{\rm M}}_{#1}(#2)}}
\newcommand{\mor}[2]{{{\rm Mor}_{#1}(\pu,#2)}}

\newcommand{\fg}{\mathfrak g}
\newcommand{\fgad}{G/P_\T}
\renewcommand{\fh}{\mathfrak h}
\newcommand{\fu}{\mathfrak u}
\newcommand{\fz}{\mathfrak z}
\newcommand{\fn}{\mathfrak n}
\newcommand{\fe}{\mathfrak e}
\newcommand{\fp}{\mathfrak p}
\newcommand{\ft}{\mathfrak t}
\newcommand{\fl}{\mathfrak l}
\newcommand{\fq}{\mathfrak q}
\newcommand{\fsl}{\mathfrak {sl}}
\newcommand{\fgl}{\mathfrak {gl}}
\newcommand{\fso}{\mathfrak {so}}
\newcommand{\fsp}{\mathfrak {sp}}
\newcommand{\ff}{\mathfrak {f}}

\newcommand{\ad}{{\rm ad}}
\newcommand{\jad}{{j^\ad}}
\newcommand{\id}{{\rm id}}


\newcommand{\dynkinadeux}[2]
{
$#1$
\setlength{\unitlength}{1.2pt}
\hspace{-3mm}
\begin{picture}(12,3)
\put(0,3){\line(1,0){10}}
\end{picture}
\hspace{-2.4mm}
$#2$
}

\newcommand{\mdynkinadeux}[2]
{
\mbox{\dynkinadeux{#1}{#2}}
}

\newcommand{\dynkingdeux}[2]
{
$#1$
\setlength{\unitlength}{1.2pt}
\hspace{-3mm}
\begin{picture}(12,3)
\put(1,.8){$<$}
\multiput(0,1.5)(0,1.5){3}{\line(1,0){10}}
\end{picture}
\hspace{-2.4mm}
$#2$
}

\newcommand{\poidsesix}[6]
{
\hspace{-.12cm}
\left (
\begin{array}{ccccc}
{} \hspace{-.2cm} #1 & {} \hspace{-.3cm} #2 & {} \hspace{-.3cm} #3 &
{} \hspace{-.3cm} #4 & {} \hspace{-.3cm} #5 \vspace{-.13cm}\\
\hspace{-.2cm} & \hspace{-.3cm} & {} \hspace{-.3cm} #6 &
{} \hspace{-.3cm} & {} \hspace{-.3cm}
\end{array}
\hspace{-.2cm}
\right )      }

\newcommand{\copoidsesix}[6]{
\hspace{-.12cm}
\left |
\begin{array}{ccccc}
{} \hspace{-.2cm} #1 & {} \hspace{-.3cm} #2 & {} \hspace{-.3cm} #3 &
{} \hspace{-.3cm} #4 & {} \hspace{-.3cm} #5 \vspace{-.13cm}\\
\hspace{-.2cm} & \hspace{-.3cm} & {} \hspace{-.3cm} #6 &
{} \hspace{-.3cm} & {} \hspace{-.3cm}
\end{array}
\hspace{-.2cm}
\right |      }

\newcommand{\poidsesept}[7]{
\hspace{-.12cm}
\left (
\begin{array}{cccccc}
{} \hspace{-.2cm} #1 & {} \hspace{-.3cm} #2 & {} \hspace{-.3cm} #3 &
{} \hspace{-.3cm} #4 & {} \hspace{-.3cm} #5 & {} \hspace{-.3cm} #6
\vspace{-.13cm}\\
\hspace{-.2cm} & \hspace{-.3cm} & {} \hspace{-.3cm} #7 &
{} \hspace{-.3cm} & {} \hspace{-.3cm}
\end{array}
\hspace{-.2cm}
\right )      }

\newcommand{\copoidsesept}[7]{
\hspace{-.12cm}
\left |
\begin{array}{cccccc}
{} \hspace{-.2cm} #1 & {} \hspace{-.3cm} #2 & {} \hspace{-.3cm} #3 &
{} \hspace{-.3cm} #4 & {} \hspace{-.3cm} #5 & {} \hspace{-.3cm} #6
\vspace{-.13cm}\\
\hspace{-.2cm} & \hspace{-.3cm} & {} \hspace{-.3cm} #7 &
{} \hspace{-.3cm} & {} \hspace{-.3cm}
\end{array}
\hspace{-.2cm}
\right |      }

\newcommand{\poidsehuit}[8]{
\hspace{-.12cm}
\left (
\begin{array}{ccccccc}
{} \hspace{-.2cm} #1 & {} \hspace{-.3cm} #2 & {} \hspace{-.3cm} #3 &
{} \hspace{-.3cm} #4 & {} \hspace{-.3cm} #5 & {} \hspace{-.3cm} #6 &
{} \hspace{-.3cm} #7   \vspace{-.13cm}\\
\hspace{-.2cm} & \hspace{-.3cm} & {} \hspace{-.3cm} #8 &
{} \hspace{-.3cm} & {} \hspace{-.3cm}
\end{array}
\hspace{-.2cm}
\right )      }

\newcommand{\copoidsehuit}[8]{
\hspace{-.12cm}
\left |
\begin{array}{cccccc}
{} \hspace{-.2cm} #1 & {} \hspace{-.3cm} #2 & {} \hspace{-.3cm} #3 &
{} \hspace{-.3cm} #4 & {} \hspace{-.3cm} #5 & {} \hspace{-.3cm} #6 &
{} \hspace{-.3cm} #7  \vspace{-.13cm}\\
\hspace{-.2cm} & \hspace{-.3cm} & {} \hspace{-.3cm} #8 &
{} \hspace{-.3cm} & {} \hspace{-.3cm}
\end{array}
\hspace{-.2cm}
\right |      }

\newcommand{\im}{\mathtt{Im}}


\def\cA{{\cal A}} \def\cC{{\cal C}} \def\cD{{\cal D}} \def\cE{{\cal E}}
\def\cF{{\cal F}} \def\cG{{\cal G}} \def\cH{{\cal H}} \def\cI{{\cal I}}
\def\cK{{\cal K}} \def\cL{{\cal L}} \def\cM{{\cal M}} \def\cN{{\cal N}}
\def\cO{{\cal O}}
\def\cP{{\cal P}} \def\cQ{{\cal Q}} \def\cT{{\cal T}} \def\cU{{\cal U}}
\def\cV{{\cal V}} \def\cX{{\cal X}} \def\cY{{\cal Y}} \def\cZ{{\cal Z}}

 \title{Computations of quantum Giambelli formulas in some exceptional
homogeneous spaces}
 \author{P.E. Chaput, N. Perrin}

\maketitle

\begin{abstract}
This paper gathers results obtained with a software written in Java,
accessible on the web-page\\
\texttt{www.math.sciences.univ-nantes.fr/}\verb|~|\texttt{chaput/quiver-demo.html}.
For exceptional minuscule, quasi-minuscule, cominuscule or adjoint
homogeneous spaces, the Schubert cells are expressed as polynomials in the
generators of the quantum cohomology algebra.
\end{abstract}

\section{Notations}

Let us explain our notations.
We consider a minuscule, quasi-minuscule, cominuscule or adjoint
rational homogeneous space $G/P$ (for the definition of these spaces,
see \cite{CP2}), and moreover we assume that $G$ is an exceptional group,
namely that its type is one of $E_6,E_7,E_8,F_4,G_2$.

Let $r$ denote the rank of $G$.
Let $(\alpha_i)_{1 \leq i \leq r}$ be the basis of the root
system of $G$ chosen in \cite{bourbaki}.
An element $\alpha = a_1 \alpha_1 + \cdots + a_r \alpha_r$ of the root lattice
will be denoted $(a_1,\ldots,a_r)$. Let $\varpi$ be the fundamental weight
corresponding to $P$.
The Schubert cells are well-known to be
parametrised by the coset $W/W_P$. We denote each Schubert cell $\sigma_w$,
with $w \in W/W_P$, with $\sigma_\alpha$, where
$\alpha = \varpi - w(\varpi)$.

We choose some generators of $H^*(G/P)$ and we give,
for each $w \in W/W_P$, an expression of $\sigma_w$ as a
polynomial in these generators in the
quantum cohomology algebra (thus the choice we made of the generators
is visible in these expressions).

\tableofcontents

\section{$E_6$}

\subsection{$E_6/P_1$}


\noindent
Schubert cells in degree $0$:

\noindent
$
\begin{array}{rcl}
\sigma_{(0,0,0,0,0,0)} & = & 1\ \\ 
\end{array}
\vspace{.5cm}
$

\noindent
Schubert cells in degree $1$:

\noindent
$
\begin{array}{rcl}
\sigma_{(1,0,0,0,0,0)} & = & h\\ 
\end{array}
\vspace{.5cm}
$

\noindent
Schubert cells in degree $2$:

\noindent
$
\begin{array}{rcl}
\sigma_{(1,0,1,0,0,0)} & = & h^{2}\\ 
\end{array}
\vspace{.5cm}
$

\noindent
Schubert cells in degree $3$:

\noindent
$
\begin{array}{rcl}
\sigma_{(1,0,1,1,0,0)} & = & h^{3}\\ 
\end{array}
\vspace{.5cm}
$

\noindent
Schubert cells in degree $4$:

\noindent
$
\begin{array}{rcl}
\sigma_{(1,0,1,1,1,0)} & = & h^{4} - s\\ 
\end{array}
\\
\begin{array}{rcl}
\sigma_{(1,1,1,1,0,0)} & = & s\\ 
\end{array}
\vspace{.5cm}
$

\noindent
Schubert cells in degree $5$:

\noindent
$
\begin{array}{rcl}
\sigma_{(1,0,1,1,1,1)} & = & h^{5} - 2\ hs\\ 
\end{array}
\\
\begin{array}{rcl}
\sigma_{(1,1,1,1,1,0)} & = & hs\\ 
\end{array}
\vspace{.5cm}
$

\noindent
Schubert cells in degree $6$:

\noindent
$
\begin{array}{rcl}
\sigma_{(1,1,1,1,1,1)} & = & h^{6} - 2\ h^{2}s\\ 
\end{array}
\\
\begin{array}{rcl}
\sigma_{(1,1,1,2,1,0)} & = & -h^{6} + 3\ h^{2}s\\ 
\end{array}
\vspace{.5cm}
$

\noindent
Schubert cells in degree $7$:

\noindent
$
\begin{array}{rcl}
\sigma_{(1,1,1,2,1,1)} & = & h^{7} - 2\ h^{3}s\\ 
\end{array}
\\
\begin{array}{rcl}
\sigma_{(1,1,2,2,1,0)} & = & -2\ h^{7} + 5\ h^{3}s\\ 
\end{array}
\vspace{.5cm}
$

\noindent
Schubert cells in degree $8$:

\noindent
$
\begin{array}{rcl}
\sigma_{(1,1,1,2,2,1)} & = & 2\ h^{8} - 5\ h^{4}s + s^{2}\\ 
\end{array}
\\
\begin{array}{rcl}
\sigma_{(1,1,2,2,1,1)} & = & -h^{8} + 3\ h^{4}s - s^{2}\\ 
\end{array}
\\
\begin{array}{rcl}
\sigma_{(2,1,2,2,1,0)} & = & -h^{8} + 2\ h^{4}s + s^{2}\\ 
\end{array}
\vspace{.5cm}
$

\noindent
Schubert cells in degree $9$:

\noindent
$
\begin{array}{rcl}
\sigma_{(1,1,2,2,2,1)} & = & 4/3\ h^{9} - 3\ h^{5}s\\ 
\end{array}
\\
\begin{array}{rcl}
\sigma_{(2,1,2,2,1,1)} & = & -5/3\ h^{9} + 4\ h^{5}s\\ 
\end{array}
\vspace{.5cm}
$

\noindent
Schubert cells in degree $10$:

\noindent
$
\begin{array}{rcl}
\sigma_{(1,1,2,3,2,1)} & = & 3\ h^{10} - 7\ h^{6}s\\ 
\end{array}
\\
\begin{array}{rcl}
\sigma_{(2,1,2,2,2,1)} & = & -5/3\ h^{10} + 4\ h^{6}s\\ 
\end{array}
\vspace{.5cm}
$

\noindent
Schubert cells in degree $11$:

\noindent
$
\begin{array}{rcl}
\sigma_{(1,2,2,3,2,1)} & = & 14/3\ h^{11} - 11\ h^{7}s\\ 
\end{array}
\\
\begin{array}{rcl}
\sigma_{(2,1,2,3,2,1)} & = & -5/3\ h^{11} + 4\ h^{7}s\\ 
\end{array}
\vspace{.5cm}
$

\noindent
Schubert cells in degree $12$:

\noindent
$
\begin{array}{rcl}
\sigma_{(2,2,2,3,2,1)} & = & -1/3\ h^{12} + h^{8}s - s^{3}\\ 
\end{array}
\\
\begin{array}{rcl}
\sigma_{(2,1,3,3,2,1)} & = & -4/3\ h^{12} + 3\ h^{8}s + s^{3}\\ 
\end{array}
\vspace{.5cm}
$

\noindent
Schubert cells in degree $13$:

\noindent
$
\begin{array}{rcl}
\sigma_{(2,2,3,3,2,1)} & = & -8/3\ h^{13} + 19/3\ h^{9}s\\ 
\end{array}
\vspace{.5cm}
$

\noindent
Schubert cells in degree $14$:

\noindent
$
\begin{array}{rcl}
\sigma_{(2,2,3,4,2,1)} & = & -19/3\ h^{14} + 15\ h^{10}s\\ 
\end{array}
\vspace{.5cm}
$

\noindent
Schubert cells in degree $15$:

\noindent
$
\begin{array}{rcl}
\sigma_{(2,2,3,4,3,1)} & = & -10\ h^{15} + 71/3\ h^{11}s\\ 
\end{array}
\vspace{.5cm}
$

\noindent
Schubert cells in degree $16$:

\noindent
$
\begin{array}{rcl}
\sigma_{(2,2,3,4,3,2)} & = & -17/3\ h^{16} + 40/3\ h^{12}s + s^{4}\\ 
\end{array}
\vspace{.5cm}
$

\subsection{$E_6/P_2$}


\noindent
Schubert cells in degree $0$:

\noindent
$
\begin{array}{rcl}
\sigma_{(0,0,0,0,0,0)} & = & 1\ \\ 
\end{array}
\vspace{.5cm}
$

\noindent
Schubert cells in degree $1$:

\noindent
$
\begin{array}{rcl}
\sigma_{(0,1,0,0,0,0)} & = & h\\ 
\end{array}
\vspace{.5cm}
$

\noindent
Schubert cells in degree $2$:

\noindent
$
\begin{array}{rcl}
\sigma_{(0,1,0,1,0,0)} & = & h^{2}\\ 
\end{array}
\vspace{.5cm}
$

\noindent
Schubert cells in degree $3$:

\noindent
$
\begin{array}{rcl}
\sigma_{(0,1,1,1,0,0)} & = & s\\ 
\end{array}
\\
\begin{array}{rcl}
\sigma_{(0,1,0,1,1,0)} & = & h^{3} - s\\ 
\end{array}
\vspace{.5cm}
$

\noindent
Schubert cells in degree $4$:

\noindent
$
\begin{array}{rcl}
\sigma_{(1,1,1,1,0,0)} & = & t\\ 
\end{array}
\\
\begin{array}{rcl}
\sigma_{(0,1,1,1,1,0)} & = & hs - t\\ 
\end{array}
\\
\begin{array}{rcl}
\sigma_{(0,1,0,1,1,1)} & = & h^{4} - 2\ hs + t\\ 
\end{array}
\vspace{.5cm}
$

\noindent
Schubert cells in degree $5$:

\noindent
$
\begin{array}{rcl}
\sigma_{(1,1,1,1,1,0)} & = & ht\\ 
\end{array}
\\
\begin{array}{rcl}
\sigma_{(0,1,1,2,1,0)} & = & -h^{5} + 3\ h^{2}s - 3\ ht\\ 
\end{array}
\\
\begin{array}{rcl}
\sigma_{(0,1,1,1,1,1)} & = & h^{5} - 2\ h^{2}s + ht\\ 
\end{array}
\vspace{.5cm}
$

\noindent
Schubert cells in degree $6$:

\noindent
$
\begin{array}{rcl}
\sigma_{(1,1,1,2,1,0)} & = & -1/2\ h^{6} + h^{3}s - 1/2\ h^{2}t + 1/2\ s^{2}\\ 
\end{array}
\\
\begin{array}{rcl}
\sigma_{(0,2,1,2,1,0)} & = & -h^{6} + 3\ h^{3}s - 2\ h^{2}t - s^{2}\\ 
\end{array}
\\
\begin{array}{rcl}
\sigma_{(1,1,1,1,1,1)} & = & 1/2\ h^{6} - h^{3}s + 3/2\ h^{2}t - 1/2\ s^{2}\\ 
\end{array}
\\
\begin{array}{rcl}
\sigma_{(0,1,1,2,1,1)} & = & 1/2\ h^{6} - h^{3}s - 1/2\ h^{2}t + 1/2\ s^{2}\\ 
\end{array}
\vspace{.5cm}
$

\noindent
Schubert cells in degree $7$:

\noindent
$
\begin{array}{rcl}
\sigma_{(1,2,1,2,1,0)} & = & -1/2\ h^{7} + h^{4}s - 1/2\ h^{3}t + 1/2\ hs^{2} - st\\ 
\end{array}
\\
\begin{array}{rcl}
\sigma_{(1,1,2,2,1,0)} & = & -1/2\ h^{7} + h^{4}s - 3/2\ h^{3}t + 1/2\ hs^{2} + st\\ 
\end{array}
\\
\begin{array}{rcl}
\sigma_{(1,1,1,2,1,1)} & = & 1/2\ h^{7} - h^{4}s + 3/2\ h^{3}t - 1/2\ hs^{2}\\ 
\end{array}
\\
\begin{array}{rcl}
\sigma_{(0,2,1,2,1,1)} & = & -1/2\ h^{7} + 2\ h^{4}s - 3/2\ h^{3}t - 3/2\ hs^{2} + st\\ 
\end{array}
\\
\begin{array}{rcl}
\sigma_{(0,1,1,2,2,1)} & = & 1/2\ h^{7} - 2\ h^{4}s - 1/2\ h^{3}t + 5/2\ hs^{2} - st\\ 
\end{array}
\vspace{.5cm}
$

\noindent
Schubert cells in degree $8$:

\noindent
$
\begin{array}{rcl}
\sigma_{(1,2,2,2,1,0)} & = & -1/3\ h^{8} - h^{4}t + 2\ h^{2}s^{2} - hst\\ 
\end{array}
\\
\begin{array}{rcl}
\sigma_{(1,2,1,2,1,1)} & = & -1/6\ h^{8} + h^{5}s + 1/2\ h^{4}t - 3/2\ h^{2}s^{2}\\ 
\end{array}
\\
\begin{array}{rcl}
\sigma_{(1,1,2,2,1,1)} & = & -1/6\ h^{8} + h^{5}s - 1/2\ h^{4}t - 3/2\ h^{2}s^{2} + 2\ hst\\ 
\end{array}
\\
\begin{array}{rcl}
\sigma_{(1,1,1,2,2,1)} & = & 5/6\ h^{8} - 3\ h^{5}s + 3/2\ h^{4}t + 5/2\ h^{2}s^{2} - 2\ hst\\ 
\end{array}
\\
\begin{array}{rcl}
\sigma_{(0,2,1,2,2,1)} & = & -1/3\ h^{8} + h^{5}s - 2\ h^{4}t + hst\\ 
\end{array}
\vspace{.5cm}
$

\noindent
Schubert cells in degree $9$:

\noindent
$
\begin{array}{rcl}
\sigma_{(1,2,2,3,1,0)} & = & 1/4\ h^{9} - 5/2\ h^{6}s - 1/4\ h^{5}t + 19/4\ h^{3}s^{2} - 3\ h^{2}st\\ 
\end{array}
\\
\begin{array}{rcl}
\sigma_{(1,2,2,2,1,1)} & = & -7/12\ h^{9} + 5/2\ h^{6}s - 3/4\ h^{5}t - 11/4\ h^{3}s^{2} + 2\ h^{2}st\\ 
\end{array}
\\
\begin{array}{rcl}
\sigma_{(1,2,1,2,2,1)} & = & 5/12\ h^{9} - 3/2\ h^{6}s + 5/4\ h^{5}t + 5/4\ h^{3}s^{2} - 2\ h^{2}st\\ 
\end{array}
\\
\begin{array}{rcl}
\sigma_{(1,1,2,2,2,1)} & = & 5/12\ h^{9} - 3/2\ h^{6}s + 1/4\ h^{5}t + 5/4\ h^{3}s^{2}\\ 
\end{array}
\\
\begin{array}{rcl}
\sigma_{(0,2,1,3,2,1)} & = & -3/4\ h^{9} + 5/2\ h^{6}s - 13/4\ h^{5}t - 5/4\ h^{3}s^{2} + 3\ h^{2}st\\ 
\end{array}
\vspace{.5cm}
$

\noindent
Schubert cells in degree $10$:

\noindent
$
\begin{array}{rcl}
\sigma_{(1,2,2,3,2,0)} & = & 7/12\ h^{10} - 7/2\ h^{7}s + 1/4\ h^{6}t + 21/4\ h^{4}s^{2} - 3\ h^{3}st - 1/2\ s^{2}t\\ 
\end{array}
\\
\begin{array}{rcl}
\sigma_{(1,2,2,3,1,1)} & = & -1/3\ h^{10} + h^{7}s - 1/2\ h^{6}t - 1/2\ h^{4}s^{2} + 1/2\ s^{2}t\\ 
\end{array}
\\
\begin{array}{rcl}
\sigma_{(1,2,2,2,2,1)} & = & -1/4\ h^{10} + 3/2\ h^{7}s - 1/4\ h^{6}t - 9/4\ h^{4}s^{2} + 2\ h^{3}st - 1/2\ s^{2}t\\ 
\end{array}
\\
\begin{array}{rcl}
\sigma_{(1,2,1,3,2,1)} & = & 2/3\ h^{10} - 3\ h^{7}s + 3/2\ h^{6}t + 7/2\ h^{4}s^{2} - 4\ h^{3}st + 1/2\ s^{2}t\\ 
\end{array}
\\
\begin{array}{rcl}
\sigma_{(0,2,2,3,2,1)} & = & -17/12\ h^{10} + 11/2\ h^{7}s - 19/4\ h^{6}t - 19/4\ h^{4}s^{2} + 7\ h^{3}st - 1/2\ s^{2}t\\ 
\end{array}
\\
\begin{array}{rcl}
\sigma_{(1,1,2,3,2,1)} & = & 2/3\ h^{10} - 3\ h^{7}s + 1/2\ h^{6}t + 7/2\ h^{4}s^{2} - 2\ h^{3}st + 1/2\ s^{2}t\\ 
\end{array}
\vspace{.5cm}
$

\noindent
Schubert cells in degree $11$:

\noindent
$
\begin{array}{rcl}
\sigma_{(1,2,2,3,2,2)} & = & -5/3\ h^{11} + 22/3\ h^{8}s - 5\ h^{7}t - 8\ h^{5}s^{2} + 12\ h^{4}st - 4\ hs^{2}t + q\\ 
\end{array}
\\
\begin{array}{rcl}
\sigma_{(1,2,2,3,3,1)} & = & 47/12\ h^{11} - 109/6\ h^{8}s + 41/4\ h^{7}t + 85/4\ h^{5}s^{2} - 27\ h^{4}st + 15/2\ hs^{2}t\\
 & &  - 2\ q\\ 
\end{array}
\\
\begin{array}{rcl}
\sigma_{(1,2,3,3,2,1)} & = & 21/4\ h^{11} - 143/6\ h^{8}s + 53/4\ h^{7}t + 109/4\ h^{5}s^{2} - 33\ h^{4}st + 15/2\ hs^{2}t\\
 & &  - 2\ q\\ 
\end{array}
\\
\begin{array}{rcl}
\sigma_{(1,2,2,4,2,1)} & = & -13/2\ h^{11} + 30\ h^{8}s - 16\ h^{7}t - 35\ h^{5}s^{2} + 42\ h^{4}st - 21/2\ hs^{2}t + 3\ q\\ 
\end{array}
\\
\begin{array}{rcl}
\sigma_{(2,2,2,3,2,1)} & = & -10/3\ h^{11} + 44/3\ h^{8}s - 9\ h^{7}t - 16\ h^{5}s^{2} + 20\ h^{4}st - 4\ hs^{2}t + q\\ 
\end{array}
\\
\begin{array}{rcl}
\sigma_{(1,3,2,3,2,1)} & = & 43/12\ h^{11} - 33/2\ h^{8}s + 33/4\ h^{7}t + 77/4\ h^{5}s^{2} - 22\ h^{4}st + 11/2\ hs^{2}t\\
 & &  - 2\ q\\ 
\end{array}
\vspace{.5cm}
$

\noindent
Schubert cells in degree $12$:

\noindent
$
\begin{array}{rcl}
\sigma_{(1,2,2,3,3,2)} & = & 17/48\ h^{12} - 15/8\ h^{9}s + 7/48\ h^{8}t + 41/16\ h^{6}s^{2} - h^{5}st - 3/4\ h^{2}s^{2}t\\ 
\end{array}
\\
\begin{array}{rcl}
\sigma_{(1,2,2,4,3,1)} & = & -23/48\ h^{12} + 17/8\ h^{9}s - 3/16\ h^{8}t - 39/16\ h^{6}s^{2} + 7/4\ h^{2}s^{2}t\\ 
\end{array}
\\
\begin{array}{rcl}
\sigma_{(1,2,3,4,2,1)} & = & 121/48\ h^{12} - 87/8\ h^{9}s + 109/16\ h^{8}t + 185/16\ h^{6}s^{2} - 14\ h^{5}st\\
 & &  + 7/4\ h^{2}s^{2}t\\ 
\end{array}
\\
\begin{array}{rcl}
\sigma_{(2,2,3,3,2,1)} & = & -21/16\ h^{12} + 131/24\ h^{9}s - 185/48\ h^{8}t - 87/16\ h^{6}s^{2} + 7\ h^{5}st\\
 & &  - 3/4\ h^{2}s^{2}t\\ 
\end{array}
\\
\begin{array}{rcl}
\sigma_{(1,3,2,4,2,1)} & = & -119/48\ h^{12} + 89/8\ h^{9}s - 115/16\ h^{8}t - 199/16\ h^{6}s^{2} + 17\ h^{5}st\\
 & &  - 17/4\ h^{2}s^{2}t\\ 
\end{array}
\vspace{.5cm}
$

\noindent
Schubert cells in degree $13$:

\noindent
$
\begin{array}{rcl}
\sigma_{(1,2,2,4,3,2)} & = & 17/48\ h^{13} - 15/8\ h^{10}s + 7/48\ h^{9}t + 41/16\ h^{7}s^{2} - h^{6}st - 3/4\ h^{3}s^{2}t\\ 
\end{array}
\\
\begin{array}{rcl}
\sigma_{(1,2,3,4,3,1)} & = & 15/4\ h^{13} - 49/3\ h^{10}s + 34/3\ h^{9}t + 35/2\ h^{7}s^{2} - 25\ h^{6}st + 25/4\ h^{3}s^{2}t\\ 
\end{array}
\\
\begin{array}{rcl}
\sigma_{(2,2,3,4,2,1)} & = & -21/16\ h^{13} + 131/24\ h^{10}s - 185/48\ h^{9}t - 87/16\ h^{7}s^{2} + 7\ h^{6}st\\
 & &  - 3/4\ h^{3}s^{2}t\\ 
\end{array}
\\
\begin{array}{rcl}
\sigma_{(1,3,2,4,3,1)} & = & -55/12\ h^{13} + 61/3\ h^{10}s - 35/3\ h^{9}t - 45/2\ h^{7}s^{2} + 26\ h^{6}st - 15/4\ h^{3}s^{2}t\\ 
\end{array}
\\
\begin{array}{rcl}
\sigma_{(1,3,3,4,2,1)} & = & 1/12\ h^{13} - 2/3\ h^{9}t - 1/2\ h^{7}s^{2} + 4\ h^{6}st - 15/4\ h^{3}s^{2}t\\ 
\end{array}
\vspace{.5cm}
$

\noindent
Schubert cells in degree $14$:

\noindent
$
\begin{array}{rcl}
\sigma_{(1,2,3,4,3,2)} & = & 305/48\ h^{14} - 677/24\ h^{11}s + 775/48\ h^{10}t + 1499/48\ h^{8}s^{2} - 36\ h^{7}st\\
 & &  + 21/4\ h^{4}s^{2}t - sq\\ 
\end{array}
\\
\begin{array}{rcl}
\sigma_{(2,2,3,4,3,1)} & = & -2\ h^{14} + 26/3\ h^{11}s - 4\ h^{10}t - 28/3\ h^{8}s^{2} + 7\ h^{7}st + 2\ h^{4}s^{2}t + sq\\ 
\end{array}
\\
\begin{array}{rcl}
\sigma_{(1,3,2,4,3,2)} & = & -6\ h^{14} + 79/3\ h^{11}s - 16\ h^{10}t - 86/3\ h^{8}s^{2} + 35\ h^{7}st - 6\ h^{4}s^{2}t + sq\\ 
\end{array}
\\
\begin{array}{rcl}
\sigma_{(1,3,3,4,3,1)} & = & -29/48\ h^{14} + 77/24\ h^{11}s - 13/16\ h^{10}t - 211/48\ h^{8}s^{2} + 4\ h^{7}st - h^{4}s^{2}t\\ 
\end{array}
\\
\begin{array}{rcl}
\sigma_{(2,3,3,4,2,1)} & = & 11/16\ h^{14} - 77/24\ h^{11}s + 7/48\ h^{10}t + 187/48\ h^{8}s^{2} - 11/4\ h^{4}s^{2}t - sq\\ 
\end{array}
\vspace{.5cm}
$

\noindent
Schubert cells in degree $15$:

\noindent
$
\begin{array}{rcl}
\sigma_{(2,2,3,4,3,2)} & = & 3\ h^{15} - 40/3\ h^{12}s + 29/3\ h^{11}t + 44/3\ h^{9}s^{2} - 70/3\ h^{8}st + 8\ h^{5}s^{2}t + tq\\ 
\end{array}
\\
\begin{array}{rcl}
\sigma_{(1,3,3,4,3,2)} & = & 49/96\ h^{15} - 109/48\ h^{12}s - 41/96\ h^{11}t + 251/96\ h^{9}s^{2} + 5/2\ h^{8}st\\
 & &  - 35/8\ h^{5}s^{2}t - tq\\ 
\end{array}
\\
\begin{array}{rcl}
\sigma_{(1,3,3,5,3,1)} & = & -173/96\ h^{15} + 139/16\ h^{12}s - 17/32\ h^{11}t - 349/32\ h^{9}s^{2} + 3/2\ h^{8}st\\
 & &  + 49/8\ h^{5}s^{2}t + 3\ tq\\ 
\end{array}
\\
\begin{array}{rcl}
\sigma_{(2,3,3,4,3,1)} & = & -69/32\ h^{15} + 451/48\ h^{12}s - 649/96\ h^{11}t - 965/96\ h^{9}s^{2} + 91/6\ h^{8}st\\
 & &  - 35/8\ h^{5}s^{2}t - tq\\ 
\end{array}
\vspace{.5cm}
$

\noindent
Schubert cells in degree $16$:

\noindent
$
\begin{array}{rcl}
\sigma_{(1,3,3,5,3,2)} & = & 3/32\ h^{16} - 1/16\ h^{13}s - 95/96\ h^{12}t - 13/32\ h^{10}s^{2} + 21/4\ h^{9}st\\
 & &  - 81/16\ h^{6}s^{2}t\\ 
\end{array}
\\
\begin{array}{rcl}
\sigma_{(2,3,3,4,3,2)} & = & 121/36\ h^{16} - 359/24\ h^{13}s + 457/48\ h^{12}t + 265/16\ h^{10}s^{2} - 89/4\ h^{9}st\\
 & &  + 89/16\ h^{6}s^{2}t\\ 
\end{array}
\\
\begin{array}{rcl}
\sigma_{(2,3,3,5,3,1)} & = & -599/96\ h^{16} + 1325/48\ h^{13}s - 1535/96\ h^{12}t - 973/32\ h^{10}s^{2} + 141/4\ h^{9}st\\
 & &  - 81/16\ h^{6}s^{2}t\\ 
\end{array}
\vspace{.5cm}
$

\noindent
Schubert cells in degree $17$:

\noindent
$
\begin{array}{rcl}
\sigma_{(1,3,3,5,4,2)} & = & -733/288\ h^{17} + 559/48\ h^{14}s - 1037/96\ h^{13}t - 1265/96\ h^{11}s^{2} + 89/3\ h^{10}st\\
 & &  - 31/2\ h^{7}s^{2}t - s^{2}q\\ 
\end{array}
\\
\begin{array}{rcl}
\sigma_{(2,3,3,5,3,2)} & = & 335/72\ h^{17} - 985/48\ h^{14}s + 1351/96\ h^{13}t + 2149/96\ h^{11}s^{2} - 781/24\ h^{10}st\\
 & &  + 295/32\ h^{7}s^{2}t + 1/2\ s^{2}q\\ 
\end{array}
\\
\begin{array}{rcl}
\sigma_{(2,3,4,5,3,1)} & = & -3613/288\ h^{17} + 885/16\ h^{14}s - 1103/32\ h^{13}t - 5809/96\ h^{11}s^{2} + 77\ h^{10}st\\
 & &  - 31/2\ h^{7}s^{2}t - s^{2}q\\ 
\end{array}
\vspace{.5cm}
$

\noindent
Schubert cells in degree $18$:

\noindent
$
\begin{array}{rcl}
\sigma_{(2,3,3,5,4,2)} & = & 4801/384\ h^{18} - 3507/64\ h^{15}s + 12971/384\ h^{14}t + 7613/128\ h^{12}s^{2}\\
 & &  - 293/4\ h^{11}st + 595/48\ h^{8}s^{2}t + stq\\ 
\end{array}
\\
\begin{array}{rcl}
\sigma_{(2,3,4,5,3,2)} & = & -10909/1152\ h^{18} + 2661/64\ h^{15}s - 9485/384\ h^{14}t - 17393/384\ h^{12}s^{2}\\
 & &  + 319/6\ h^{11}st - 59/8\ h^{8}s^{2}t - stq\\ 
\end{array}
\vspace{.5cm}
$

\noindent
Schubert cells in degree $19$:

\noindent
$
\begin{array}{rcl}
\sigma_{(2,3,4,5,4,2)} & = & 1429/288\ h^{19} - 1567/72\ h^{16}s + 2099/144\ h^{15}t + 283/12\ h^{13}s^{2}\\
 & &  - 787/24\ h^{12}st + 787/96\ h^{9}s^{2}t\\ 
\end{array}
\vspace{.5cm}
$

\noindent
Schubert cells in degree $20$:

\noindent
$
\begin{array}{rcl}
\sigma_{(2,3,4,6,4,2)} & = & 27497/2304\ h^{20} - 6721/128\ h^{17}s + 8899/256\ h^{16}t + 43901/768\ h^{14}s^{2}\\
 & &  - 473/6\ h^{13}st + 473/24\ h^{10}s^{2}t\\ 
\end{array}
\vspace{.5cm}
$

\noindent
Schubert cells in degree $21$:

\noindent
$
\begin{array}{rcl}
\sigma_{(2,4,4,6,4,2)} & = & 34663/2304\ h^{21} - 76567/1152\ h^{18}s + 105349/2304\ h^{17}t + 55811/768\ h^{15}s^{2}\\
 & &  - 427/4\ h^{14}st + 1505/48\ h^{11}s^{2}t - 1/2\ s^{2}tq\\ 
\end{array}
\vspace{.5cm}
$

\newpage

\section{$E_7$}

\subsection{$E_7/P_1$}


\noindent
Schubert cells in degree $0$:

\noindent
$
\begin{array}{rcl}
\sigma_{(0,0,0,0,0,0,0)} & = & 1\ \\ 
\end{array}
\vspace{.5cm}
$

\noindent
Schubert cells in degree $1$:

\noindent
$
\begin{array}{rcl}
\sigma_{(1,0,0,0,0,0,0)} & = & h\\ 
\end{array}
\vspace{.5cm}
$

\noindent
Schubert cells in degree $2$:

\noindent
$
\begin{array}{rcl}
\sigma_{(1,1,0,0,0,0,0)} & = & h^{2}\\ 
\end{array}
\vspace{.5cm}
$

\noindent
Schubert cells in degree $3$:

\noindent
$
\begin{array}{rcl}
\sigma_{(1,1,1,0,0,0,0)} & = & h^{3}\\ 
\end{array}
\vspace{.5cm}
$

\noindent
Schubert cells in degree $4$:

\noindent
$
\begin{array}{rcl}
\sigma_{(1,1,1,1,0,0,0)} & = & s\\ 
\end{array}
\\
\begin{array}{rcl}
\sigma_{(1,1,1,0,1,0,0)} & = & h^{4} - s\\ 
\end{array}
\vspace{.5cm}
$

\noindent
Schubert cells in degree $5$:

\noindent
$
\begin{array}{rcl}
\sigma_{(1,1,1,1,1,0,0)} & = & hs\\ 
\end{array}
\\
\begin{array}{rcl}
\sigma_{(1,1,1,0,1,1,0)} & = & h^{5} - 2\ hs\\ 
\end{array}
\vspace{.5cm}
$

\noindent
Schubert cells in degree $6$:

\noindent
$
\begin{array}{rcl}
\sigma_{(1,1,2,1,1,0,0)} & = & -h^{6} + 3\ h^{2}s + t\\ 
\end{array}
\\
\begin{array}{rcl}
\sigma_{(1,1,1,1,1,1,0)} & = & h^{6} - 2\ h^{2}s - t\\ 
\end{array}
\\
\begin{array}{rcl}
\sigma_{(1,1,1,0,1,1,1)} & = & t\\ 
\end{array}
\vspace{.5cm}
$

\noindent
Schubert cells in degree $7$:

\noindent
$
\begin{array}{rcl}
\sigma_{(1,2,2,1,1,0,0)} & = & -2\ h^{7} + 5\ h^{3}s + 3\ ht\\ 
\end{array}
\\
\begin{array}{rcl}
\sigma_{(1,1,2,1,1,1,0)} & = & h^{7} - 2\ h^{3}s - 2\ ht\\ 
\end{array}
\\
\begin{array}{rcl}
\sigma_{(1,1,1,1,1,1,1)} & = & ht\\ 
\end{array}
\vspace{.5cm}
$

\noindent
Schubert cells in degree $8$:

\noindent
$
\begin{array}{rcl}
\sigma_{(2,2,2,1,1,0,0)} & = & -h^{8} + 2\ h^{4}s + 2\ h^{2}t + s^{2}\\ 
\end{array}
\\
\begin{array}{rcl}
\sigma_{(1,2,2,1,1,1,0)} & = & -h^{8} + 3\ h^{4}s + h^{2}t - s^{2}\\ 
\end{array}
\\
\begin{array}{rcl}
\sigma_{(1,1,2,1,2,1,0)} & = & 2\ h^{8} - 5\ h^{4}s - 4\ h^{2}t + s^{2}\\ 
\end{array}
\\
\begin{array}{rcl}
\sigma_{(1,1,2,1,1,1,1)} & = & h^{2}t\\ 
\end{array}
\vspace{.5cm}
$

\noindent
Schubert cells in degree $9$:

\noindent
$
\begin{array}{rcl}
\sigma_{(2,2,2,1,1,1,0)} & = & -h^{9} + 2\ h^{5}s + 2\ h^{3}t + hs^{2}\\ 
\end{array}
\\
\begin{array}{rcl}
\sigma_{(1,2,2,1,2,1,0)} & = & h^{9} - 2\ h^{5}s - 3\ h^{3}t - 1/2\ hs^{2}\\ 
\end{array}
\\
\begin{array}{rcl}
\sigma_{(1,2,2,1,1,1,1)} & = & -h^{9} + 3\ h^{5}s + 2\ h^{3}t - 3/2\ hs^{2}\\ 
\end{array}
\\
\begin{array}{rcl}
\sigma_{(1,1,2,1,2,1,1)} & = & h^{9} - 3\ h^{5}s - h^{3}t + 3/2\ hs^{2}\\ 
\end{array}
\vspace{.5cm}
$

\noindent
Schubert cells in degree $10$:

\noindent
$
\begin{array}{rcl}
\sigma_{(2,2,2,1,2,1,0)} & = & -h^{6}s + 5/2\ h^{2}s^{2} + st\\ 
\end{array}
\\
\begin{array}{rcl}
\sigma_{(1,2,3,1,2,1,0)} & = & h^{10} - h^{6}s - 3\ h^{4}t - 3\ h^{2}s^{2} - 2\ st\\ 
\end{array}
\\
\begin{array}{rcl}
\sigma_{(2,2,2,1,1,1,1)} & = & -h^{10} + 3\ h^{6}s + 2\ h^{4}t - 3/2\ h^{2}s^{2} - st\\ 
\end{array}
\\
\begin{array}{rcl}
\sigma_{(1,2,2,1,2,1,1)} & = & st\\ 
\end{array}
\\
\begin{array}{rcl}
\sigma_{(1,1,2,1,2,2,1)} & = & h^{10} - 3\ h^{6}s - h^{4}t + 3/2\ h^{2}s^{2} - st\\ 
\end{array}
\vspace{.5cm}
$

\noindent
Schubert cells in degree $11$:

\noindent
$
\begin{array}{rcl}
\sigma_{(2,2,3,1,2,1,0)} & = & h^{11} - 4\ h^{7}s - 2\ h^{5}t + 4\ h^{3}s^{2} + 2\ hst\\ 
\end{array}
\\
\begin{array}{rcl}
\sigma_{(1,2,3,2,2,1,0)} & = & 3\ h^{7}s - 7\ h^{3}s^{2} - 7\ hst\\ 
\end{array}
\\
\begin{array}{rcl}
\sigma_{(2,2,2,1,2,1,1)} & = & -h^{11} + 3\ h^{7}s + 2\ h^{5}t - 3/2\ h^{3}s^{2} - hst\\ 
\end{array}
\\
\begin{array}{rcl}
\sigma_{(1,2,3,1,2,1,1)} & = & -h^{5}t + 3\ hst\\ 
\end{array}
\\
\begin{array}{rcl}
\sigma_{(1,2,2,1,2,2,1)} & = & h^{11} - 3\ h^{7}s - h^{5}t + 3/2\ h^{3}s^{2} - hst\\ 
\end{array}
\vspace{.5cm}
$

\noindent
Schubert cells in degree $12$:

\noindent
$
\begin{array}{rcl}
\sigma_{(2,3,3,1,2,1,0)} & = & h^{12} - 4\ h^{8}s - 2\ h^{6}t + 7/2\ h^{4}s^{2} + 3\ h^{2}st + s^{3}\\ 
\end{array}
\\
\begin{array}{rcl}
\sigma_{(2,2,3,2,2,1,0)} & = & 2/3\ h^{12} - h^{8}s - 2/3\ h^{6}t - h^{4}s^{2} - 4\ h^{2}st - 2/3\ s^{3}\\ 
\end{array}
\\
\begin{array}{rcl}
\sigma_{(2,2,3,1,2,1,1)} & = & -2/3\ h^{12} + h^{8}s + 2/3\ h^{6}t + 3/2\ h^{4}s^{2} + 3\ h^{2}st - 1/3\ s^{3}\\ 
\end{array}
\\
\begin{array}{rcl}
\sigma_{(1,2,3,2,2,1,1)} & = & -2/3\ h^{12} + 4\ h^{8}s + 2/3\ h^{6}t - 6\ h^{4}s^{2} - 3\ h^{2}st + 2/3\ s^{3}\\ 
\end{array}
\\
\begin{array}{rcl}
\sigma_{(2,2,2,1,2,2,1)} & = & -1/3\ h^{12} + 2\ h^{8}s + 4/3\ h^{6}t - 3\ h^{4}s^{2} - 4\ h^{2}st + 1/3\ s^{3}\\ 
\end{array}
\\
\begin{array}{rcl}
\sigma_{(1,2,3,1,2,2,1)} & = & 4/3\ h^{12} - 5\ h^{8}s - 7/3\ h^{6}t + 9/2\ h^{4}s^{2} + 3\ h^{2}st - 1/3\ s^{3}\\ 
\end{array}
\vspace{.5cm}
$

\noindent
Schubert cells in degree $13$:

\noindent
$
\begin{array}{rcl}
\sigma_{(2,3,3,2,2,1,0)} & = & h^{13} - 2\ h^{9}s - h^{7}t - h^{5}s^{2} - 4\ h^{3}st + 1/2\ hs^{3}\\ 
\end{array}
\\
\begin{array}{rcl}
\sigma_{(2,3,3,1,2,1,1)} & = & -2\ h^{9}s - h^{7}t + 9/2\ h^{5}s^{2} + 7\ h^{3}st + 1/2\ hs^{3}\\ 
\end{array}
\\
\begin{array}{rcl}
\sigma_{(2,2,3,2,2,1,1)} & = & -1/3\ h^{13} + h^{9}s + 1/3\ h^{7}t - 7/6\ hs^{3}\\ 
\end{array}
\\
\begin{array}{rcl}
\sigma_{(2,2,3,1,2,2,1)} & = & -1/3\ h^{13} + 2\ h^{9}s + 4/3\ h^{7}t - 3\ h^{5}s^{2} - 4\ h^{3}st + 1/3\ hs^{3}\\ 
\end{array}
\\
\begin{array}{rcl}
\sigma_{(1,2,3,2,2,2,1)} & = & -1/3\ h^{13} + 3\ h^{9}s + 1/3\ h^{7}t - 6\ h^{5}s^{2} - 3\ h^{3}st + 11/6\ hs^{3}\\ 
\end{array}
\\
\begin{array}{rcl}
\sigma_{(1,2,3,1,3,2,1)} & = & 2\ h^{13} - 10\ h^{9}s - 4\ h^{7}t + 27/2\ h^{5}s^{2} + 10\ h^{3}st - 5/2\ hs^{3}\\ 
\end{array}
\vspace{.5cm}
$

\noindent
Schubert cells in degree $14$:

\noindent
$
\begin{array}{rcl}
\sigma_{(2,3,4,2,2,1,0)} & = & 1/3\ h^{14} + 3\ h^{10}s + h^{8}t - 10\ h^{6}s^{2} - 12\ h^{4}st + 8/3\ h^{2}s^{3}\\ 
\end{array}
\\
\begin{array}{rcl}
\sigma_{(2,3,3,2,2,1,1)} & = & 2/3\ h^{14} - 5\ h^{10}s - 2\ h^{8}t + 9\ h^{6}s^{2} + 8\ h^{4}st - 13/6\ h^{2}s^{3}\\ 
\end{array}
\\
\begin{array}{rcl}
\sigma_{(2,3,3,1,2,2,1)} & = & -2/3\ h^{14} + 3\ h^{10}s + h^{8}t - 9/2\ h^{6}s^{2} - h^{4}st + 8/3\ h^{2}s^{3}\\ 
\end{array}
\\
\begin{array}{rcl}
\sigma_{(2,2,3,2,2,2,1)} & = & -h^{14} + 6\ h^{10}s + 7/3\ h^{8}t - 9\ h^{6}s^{2} - 8\ h^{4}st + h^{2}s^{3}\\ 
\end{array}
\\
\begin{array}{rcl}
\sigma_{(2,2,3,1,3,2,1)} & = & 4/3\ h^{14} - 7\ h^{10}s - 2\ h^{8}t + 21/2\ h^{6}s^{2} + 5\ h^{4}st - 10/3\ h^{2}s^{3}\\ 
\end{array}
\\
\begin{array}{rcl}
\sigma_{(1,2,3,2,3,2,1)} & = & 2/3\ h^{14} - 3\ h^{10}s - 2\ h^{8}t + 3\ h^{6}s^{2} + 5\ h^{4}st + 5/6\ h^{2}s^{3}\\ 
\end{array}
\vspace{.5cm}
$

\noindent
Schubert cells in degree $15$:

\noindent
$
\begin{array}{rcl}
\sigma_{(2,3,4,2,3,1,0)} & = & -11/6\ h^{15} + 16\ h^{11}s + 17/3\ h^{9}t - 31\ h^{7}s^{2} - 27\ h^{5}st + 25/3\ h^{3}s^{3}\\ 
\end{array}
\\
\begin{array}{rcl}
\sigma_{(2,3,4,2,2,1,1)} & = & 13/6\ h^{15} - 13\ h^{11}s - 14/3\ h^{9}t + 21\ h^{7}s^{2} + 15\ h^{5}st - 17/3\ h^{3}s^{3}\\ 
\end{array}
\\
\begin{array}{rcl}
\sigma_{(2,3,3,2,2,2,1)} & = & -3/2\ h^{15} + 8\ h^{11}s + 8/3\ h^{9}t - 12\ h^{7}s^{2} - 7\ h^{5}st + 7/2\ h^{3}s^{3}\\ 
\end{array}
\\
\begin{array}{rcl}
\sigma_{(2,3,3,1,3,2,1)} & = & 5/6\ h^{15} - 5\ h^{11}s - 5/3\ h^{9}t + 15/2\ h^{7}s^{2} + 6\ h^{5}st - 5/6\ h^{3}s^{3}\\ 
\end{array}
\\
\begin{array}{rcl}
\sigma_{(2,2,3,2,3,2,1)} & = & 1/2\ h^{15} - 2\ h^{11}s - 1/3\ h^{9}t + 3\ h^{7}s^{2} - h^{5}st - 5/2\ h^{3}s^{3}\\ 
\end{array}
\\
\begin{array}{rcl}
\sigma_{(1,2,4,2,3,2,1)} & = & 1/6\ h^{15} - h^{11}s - 5/3\ h^{9}t + 6\ h^{5}st + 10/3\ h^{3}s^{3}\\ 
\end{array}
\vspace{.5cm}
$

\noindent
Schubert cells in degree $16$:

\noindent
$
\begin{array}{rcl}
\sigma_{(2,3,4,2,3,2,0)} & = & -4\ h^{16} + 83/3\ h^{12}s + 29/3\ h^{10}t - 97/2\ h^{8}s^{2} - 113/3\ h^{6}st + 13\ h^{4}s^{3} + 1/3\ s^{4}\\ 
\end{array}
\\
\begin{array}{rcl}
\sigma_{(2,3,4,2,3,1,1)} & = & 13/6\ h^{16} - 35/3\ h^{12}s - 4\ h^{10}t + 35/2\ h^{8}s^{2} + 32/3\ h^{6}st - 14/3\ h^{4}s^{3} - 1/3\ s^{4}\\ 
\end{array}
\\
\begin{array}{rcl}
\sigma_{(2,3,4,2,2,2,1)} & = & -4/3\ h^{12}s - 2/3\ h^{10}t + 7/2\ h^{8}s^{2} + 13/3\ h^{6}st - h^{4}s^{3} + 1/3\ s^{4}\\ 
\end{array}
\\
\begin{array}{rcl}
\sigma_{(2,3,3,2,3,2,1)} & = & -3/2\ h^{16} + 28/3\ h^{12}s + 10/3\ h^{10}t - 31/2\ h^{8}s^{2} - 34/3\ h^{6}st + 9/2\ h^{4}s^{3}\\
 & &  - 1/3\ s^{4}\\ 
\end{array}
\\
\begin{array}{rcl}
\sigma_{(2,3,4,1,3,2,1)} & = & 7/3\ h^{16} - 43/3\ h^{12}s - 5\ h^{10}t + 23\ h^{8}s^{2} + 52/3\ h^{6}st - 16/3\ h^{4}s^{3} + 1/3\ s^{4}\\ 
\end{array}
\\
\begin{array}{rcl}
\sigma_{(2,2,4,2,3,2,1)} & = & 2\ h^{16} - 34/3\ h^{12}s - 11/3\ h^{10}t + 37/2\ h^{8}s^{2} + 31/3\ h^{6}st - 7\ h^{4}s^{3} + 1/3\ s^{4}\\ 
\end{array}
\\
\begin{array}{rcl}
\sigma_{(1,3,4,2,3,2,1)} & = & -11/6\ h^{16} + 31/3\ h^{12}s + 2\ h^{10}t - 37/2\ h^{8}s^{2} - 13/3\ h^{6}st + 31/3\ h^{4}s^{3} - 1/3\ s^{4}\\ 
\end{array}
\vspace{.5cm}
$

\noindent
Schubert cells in degree $17$:

\noindent
$
\begin{array}{rcl}
\sigma_{(2,3,4,2,3,2,2)} & = & 7/2\ h^{17} - 21\ h^{13}s - 6\ h^{11}t + 75/2\ h^{9}s^{2} + 18\ h^{7}st - 39/2\ h^{5}s^{3} + 9/2\ hs^{4}\\
 & &  + q\\ 
\end{array}
\\
\begin{array}{rcl}
\sigma_{(2,3,4,2,3,3,1)} & = & -11\ h^{17} + 209/3\ h^{13}s + 65/3\ h^{11}t - 247/2\ h^{9}s^{2} - 221/3\ h^{7}st + 52\ h^{5}s^{3}\\
 & &  - 26/3\ hs^{4} - 2\ q\\ 
\end{array}
\\
\begin{array}{rcl}
\sigma_{(2,3,4,2,4,2,1)} & = & 62/3\ h^{17} - 130\ h^{13}s - 124/3\ h^{11}t + 227\ h^{9}s^{2} + 140\ h^{7}st - 535/6\ h^{5}s^{3}\\
 & &  + 25/2\ hs^{4} + 3\ q\\ 
\end{array}
\\
\begin{array}{rcl}
\sigma_{(2,3,5,2,3,2,1)} & = & -91/3\ h^{17} + 189\ h^{13}s + 181/3\ h^{11}t - 327\ h^{9}s^{2} - 202\ h^{7}st + 376/3\ h^{5}s^{3}\\
 & &  - 16\ hs^{4} - 4\ q\\ 
\end{array}
\\
\begin{array}{rcl}
\sigma_{(2,3,4,3,3,2,1)} & = & 49/3\ h^{17} - 305/3\ h^{13}s - 98/3\ h^{11}t + 175\ h^{9}s^{2} + 329/3\ h^{7}st - 196/3\ h^{5}s^{3}\\
 & &  + 49/6\ hs^{4} + 2\ q\\ 
\end{array}
\\
\begin{array}{rcl}
\sigma_{(2,4,4,2,3,2,1)} & = & 133/6\ h^{17} - 137\ h^{13}s - 130/3\ h^{11}t + 473/2\ h^{9}s^{2} + 143\ h^{7}st - 275/3\ h^{5}s^{3}\\
 & &  + 11\ hs^{4} + 3\ q\\ 
\end{array}
\\
\begin{array}{rcl}
\sigma_{(3,3,4,2,3,2,1)} & = & -12\ h^{17} + 221/3\ h^{13}s + 68/3\ h^{11}t - 255/2\ h^{9}s^{2} - 221/3\ h^{7}st + 51\ h^{5}s^{3}\\
 & &  - 17/3\ hs^{4} - 2\ q\\ 
\end{array}
\vspace{.5cm}
$

\noindent
Schubert cells in degree $18$:

\noindent
$
\begin{array}{rcl}
\sigma_{(2,3,4,2,3,3,2)} & = & -53/18\ h^{18} + 172/9\ h^{14}s + 20/3\ h^{12}t - 193/6\ h^{10}s^{2} - 218/9\ h^{8}st\\
 & &  + 143/18\ h^{6}s^{3} + 29/36\ h^{2}s^{4}\\ 
\end{array}
\\
\begin{array}{rcl}
\sigma_{(2,3,4,2,4,3,1)} & = & 29/6\ h^{18} - 89/3\ h^{14}s - 31/3\ h^{12}t + 48\ h^{10}s^{2} + 35\ h^{8}st - 65/6\ h^{6}s^{3}\\
 & &  - 25/12\ h^{2}s^{4}\\ 
\end{array}
\\
\begin{array}{rcl}
\sigma_{(2,3,5,2,4,2,1)} & = & -7/2\ h^{18} + 20\ h^{14}s + 7\ h^{12}t - 30\ h^{10}s^{2} - 65/3\ h^{8}st + 4\ h^{6}s^{3} + 7/2\ h^{2}s^{4}\\ 
\end{array}
\\
\begin{array}{rcl}
\sigma_{(2,3,5,3,3,2,1)} & = & 31/9\ h^{18} - 193/9\ h^{14}s - 22/3\ h^{12}t + 107/3\ h^{10}s^{2} + 227/9\ h^{8}st - 94/9\ h^{6}s^{3}\\
 & &  + 7/9\ h^{2}s^{4}\\ 
\end{array}
\\
\begin{array}{rcl}
\sigma_{(2,4,5,2,3,2,1)} & = & -9/2\ h^{18} + 30\ h^{14}s + 10\ h^{12}t - 54\ h^{10}s^{2} - 110/3\ h^{8}st + 22\ h^{6}s^{3} - 11/2\ h^{2}s^{4}\\ 
\end{array}
\\
\begin{array}{rcl}
\sigma_{(3,4,4,2,3,2,1)} & = & 22/3\ h^{18} - 140/3\ h^{14}s - 46/3\ h^{12}t + 163/2\ h^{10}s^{2} + 53\ h^{8}st - 94/3\ h^{6}s^{3}\\
 & &  + 65/12\ h^{2}s^{4}\\ 
\end{array}
\vspace{.5cm}
$

\noindent
Schubert cells in degree $19$:

\noindent
$
\begin{array}{rcl}
\sigma_{(2,3,4,2,4,3,2)} & = & -53/18\ h^{19} + 172/9\ h^{15}s + 20/3\ h^{13}t - 193/6\ h^{11}s^{2} - 218/9\ h^{9}st\\
 & &  + 143/18\ h^{7}s^{3} + 29/36\ h^{3}s^{4}\\ 
\end{array}
\\
\begin{array}{rcl}
\sigma_{(2,3,5,3,4,2,1)} & = & 47/9\ h^{19} - 625/18\ h^{15}s - 32/3\ h^{13}t + 196/3\ h^{11}s^{2} + 343/9\ h^{9}st\\
 & &  - 308/9\ h^{7}s^{3} + 98/9\ h^{3}s^{4}\\ 
\end{array}
\\
\begin{array}{rcl}
\sigma_{(2,3,5,2,4,3,1)} & = & 70/9\ h^{19} - 439/9\ h^{15}s - 17\ h^{13}t + 481/6\ h^{11}s^{2} + 533/9\ h^{9}st - 169/9\ h^{7}s^{3}\\
 & &  - 26/9\ h^{3}s^{4}\\ 
\end{array}
\\
\begin{array}{rcl}
\sigma_{(2,4,5,2,4,2,1)} & = & -33/2\ h^{19} + 207/2\ h^{15}s + 104/3\ h^{13}t - 351/2\ h^{11}s^{2} - 119\ h^{9}st + 57\ h^{7}s^{3}\\
 & &  - 9/2\ h^{3}s^{4}\\ 
\end{array}
\\
\begin{array}{rcl}
\sigma_{(2,4,5,3,3,2,1)} & = & -16/9\ h^{19} + 239/18\ h^{15}s + 10/3\ h^{13}t - 89/3\ h^{11}s^{2} - 116/9\ h^{9}st\\
 & &  + 214/9\ h^{7}s^{3} - 91/9\ h^{3}s^{4}\\ 
\end{array}
\\
\begin{array}{rcl}
\sigma_{(3,4,5,2,3,2,1)} & = & 124/9\ h^{19} - 781/9\ h^{15}s - 28\ h^{13}t + 907/6\ h^{11}s^{2} + 857/9\ h^{9}st - 529/9\ h^{7}s^{3}\\
 & &  + 82/9\ h^{3}s^{4}\\ 
\end{array}
\vspace{.5cm}
$

\noindent
Schubert cells in degree $20$:

\noindent
$
\begin{array}{rcl}
\sigma_{(2,3,5,3,4,3,1)} & = & 179/9\ h^{20} - 127\ h^{16}s - 371/9\ h^{14}t + 223\ h^{12}s^{2} + 143\ h^{10}st - 1573/18\ h^{8}s^{3}\\
 & &  + 91/6\ h^{4}s^{4}\\ 
\end{array}
\\
\begin{array}{rcl}
\sigma_{(2,3,5,2,4,3,2)} & = & -53/18\ h^{20} + 172/9\ h^{16}s + 20/3\ h^{14}t - 193/6\ h^{12}s^{2} - 218/9\ h^{10}st\\
 & &  + 143/18\ h^{8}s^{3} + 29/36\ h^{4}s^{4}\\ 
\end{array}
\\
\begin{array}{rcl}
\sigma_{(2,4,5,3,4,2,1)} & = & -44/3\ h^{20} + 1661/18\ h^{16}s + 275/9\ h^{14}t - 473/3\ h^{12}s^{2} - 944/9\ h^{10}st\\
 & &  + 319/6\ h^{8}s^{3} - 77/18\ h^{4}s^{4}\\ 
\end{array}
\\
\begin{array}{rcl}
\sigma_{(2,4,5,2,4,3,1)} & = & -55/6\ h^{20} + 532/9\ h^{16}s + 158/9\ h^{14}t - 332/3\ h^{12}s^{2} - 536/9\ h^{10}st\\
 & &  + 182/3\ h^{8}s^{3} - 679/36\ h^{4}s^{4}\\ 
\end{array}
\\
\begin{array}{rcl}
\sigma_{(3,4,5,2,4,2,1)} & = & 22/3\ h^{20} - 431/9\ h^{16}s - 121/9\ h^{14}t + 557/6\ h^{12}s^{2} + 409/9\ h^{10}st\\
 & &  - 341/6\ h^{8}s^{3} + 671/36\ h^{4}s^{4}\\ 
\end{array}
\\
\begin{array}{rcl}
\sigma_{(3,4,5,3,3,2,1)} & = & 116/9\ h^{20} - 79\ h^{16}s - 245/9\ h^{14}t + 128\ h^{12}s^{2} + 92\ h^{10}st - 529/18\ h^{8}s^{3}\\
 & &  - 35/6\ h^{4}s^{4}\\ 
\end{array}
\vspace{.5cm}
$

\noindent
Schubert cells in degree $21$:

\noindent
$
\begin{array}{rcl}
\sigma_{(2,3,5,3,4,3,2)} & = & 27/2\ h^{21} - 88\ h^{17}s - 51/2\ h^{15}t + 168\ h^{13}s^{2} + 87\ h^{11}st - 195/2\ h^{9}s^{3}\\
 & &  + 63/2\ h^{5}s^{4} + sq\\ 
\end{array}
\\
\begin{array}{rcl}
\sigma_{(2,4,5,3,4,3,1)} & = & 115/18\ h^{21} - 39\ h^{17}s - 283/18\ h^{15}t + 55\ h^{13}s^{2} + 56\ h^{11}st + 91/9\ h^{9}s^{3}\\
 & &  - 49/3\ h^{5}s^{4} - sq\\ 
\end{array}
\\
\begin{array}{rcl}
\sigma_{(2,4,6,3,4,2,1)} & = & -611/18\ h^{21} + 3785/18\ h^{17}s + 147/2\ h^{15}t - 1022/3\ h^{13}s^{2} - 2276/9\ h^{11}st\\
 & &  + 652/9\ h^{9}s^{3} + 161/9\ h^{5}s^{4} + 2\ sq\\ 
\end{array}
\\
\begin{array}{rcl}
\sigma_{(2,4,5,2,4,3,2)} & = & -148/9\ h^{21} + 964/9\ h^{17}s + 193/6\ h^{15}t - 1201/6\ h^{13}s^{2} - 1001/9\ h^{11}st\\
 & &  + 949/9\ h^{9}s^{3} - 1105/36\ h^{5}s^{4} - sq\\ 
\end{array}
\\
\begin{array}{rcl}
\sigma_{(3,4,5,3,4,2,1)} & = & 116/9\ h^{21} - 79\ h^{17}s - 245/9\ h^{15}t + 128\ h^{13}s^{2} + 92\ h^{11}st - 529/18\ h^{9}s^{3}\\
 & &  - 35/6\ h^{5}s^{4} - sq\\ 
\end{array}
\\
\begin{array}{rcl}
\sigma_{(3,4,5,2,4,3,1)} & = & 8/9\ h^{21} - 9\ h^{17}s + 10/9\ h^{15}t + 69/2\ h^{13}s^{2} - 13/3\ h^{11}st - 494/9\ h^{9}s^{3}\\
 & &  + 169/6\ h^{5}s^{4} + 2\ sq\\ 
\end{array}
\vspace{.5cm}
$

\noindent
Schubert cells in degree $22$:

\noindent
$
\begin{array}{rcl}
\sigma_{(2,4,6,3,4,3,1)} & = & -1561/54\ h^{22} + 4939/27\ h^{18}s + 521/9\ h^{16}t - 2912/9\ h^{14}s^{2} - 5312/27\ h^{12}st\\
 & &  + 3692/27\ h^{10}s^{3} - 715/27\ h^{6}s^{4}\\ 
\end{array}
\\
\begin{array}{rcl}
\sigma_{(2,4,5,3,4,3,2)} & = & 1249/108\ h^{22} - 8015/108\ h^{18}s - 823/36\ h^{16}t + 2435/18\ h^{14}s^{2}\\
 & &  + 4223/54\ h^{12}st - 1768/27\ h^{10}s^{3} + 461/27\ h^{6}s^{4}\\ 
\end{array}
\\
\begin{array}{rcl}
\sigma_{(3,4,5,3,4,3,1)} & = & 77/3\ h^{22} - 323/2\ h^{18}s - 961/18\ h^{16}t + 276\ h^{14}s^{2} + 550/3\ h^{12}st\\
 & &  - 559/6\ h^{10}s^{3} + 15/2\ h^{6}s^{4}\\ 
\end{array}
\\
\begin{array}{rcl}
\sigma_{(3,4,6,3,4,2,1)} & = & -481/54\ h^{22} + 1483/27\ h^{18}s + 188/9\ h^{16}t - 743/9\ h^{14}s^{2} - 1991/27\ h^{12}st\\
 & &  - 7/27\ h^{10}s^{3} + 419/27\ h^{6}s^{4}\\ 
\end{array}
\\
\begin{array}{rcl}
\sigma_{(3,4,5,2,4,3,2)} & = & -704/27\ h^{22} + 9047/54\ h^{18}s + 943/18\ h^{16}t - 5449/18\ h^{14}s^{2} - 4877/27\ h^{12}st\\
 & &  + 7501/54\ h^{10}s^{3} - 3601/108\ h^{6}s^{4}\\ 
\end{array}
\vspace{.5cm}
$

\noindent
Schubert cells in degree $23$:

\noindent
$
\begin{array}{rcl}
\sigma_{(2,4,6,3,5,3,1)} & = & -4351/36\ h^{23} + 3065/4\ h^{19}s + 8867/36\ h^{17}t - 2689/2\ h^{15}s^{2} - 5065/6\ h^{13}st\\
 & &  + 9659/18\ h^{11}s^{3} - 1111/12\ h^{7}s^{4} - 3\ tq\\ 
\end{array}
\\
\begin{array}{rcl}
\sigma_{(2,4,6,3,4,3,2)} & = & 1355/36\ h^{23} - 967/4\ h^{19}s - 301/4\ h^{17}t + 438\ h^{15}s^{2} + 1553/6\ h^{13}st\\
 & &  - 3679/18\ h^{11}s^{3} + 605/12\ h^{7}s^{4} + tq\\ 
\end{array}
\\
\begin{array}{rcl}
\sigma_{(3,4,6,3,4,3,1)} & = & 2933/54\ h^{23} - 18445/54\ h^{19}s - 679/6\ h^{17}t + 10493/18\ h^{15}s^{2}\\
 & &  + 10492/27\ h^{13}st - 5278/27\ h^{11}s^{3} + 847/54\ h^{7}s^{4} + 2\ tq\\ 
\end{array}
\\
\begin{array}{rcl}
\sigma_{(3,4,5,3,4,3,2)} & = & -704/27\ h^{23} + 9047/54\ h^{19}s + 943/18\ h^{17}t - 5449/18\ h^{15}s^{2} - 4877/27\ h^{13}st\\
 & &  + 7501/54\ h^{11}s^{3} - 3601/108\ h^{7}s^{4} - tq\\ 
\end{array}
\\
\begin{array}{rcl}
\sigma_{(3,5,6,3,4,2,1)} & = & -2299/36\ h^{23} + 1581/4\ h^{19}s + 4997/36\ h^{17}t - 637\ h^{15}s^{2} - 2869/6\ h^{13}st\\
 & &  + 2279/18\ h^{11}s^{3} + 473/12\ h^{7}s^{4} - 3\ tq\\ 
\end{array}
\vspace{.5cm}
$

\noindent
Schubert cells in degree $24$:

\noindent
$
\begin{array}{rcl}
\sigma_{(2,4,6,3,5,3,2)} & = & 143/18\ h^{24} - 176/3\ h^{20}s - 34/3\ h^{18}t + 571/4\ h^{16}s^{2} + 361/9\ h^{14}st\\
 & &  - 1291/9\ h^{12}s^{3} + 271/4\ h^{8}s^{4}\\ 
\end{array}
\\
\begin{array}{rcl}
\sigma_{(3,4,6,3,5,3,1)} & = & -1972/27\ h^{24} + 12806/27\ h^{20}s + 424/3\ h^{18}t - 16007/18\ h^{16}s^{2}\\
 & &  - 13093/27\ h^{14}st + 25931/54\ h^{12}s^{3} - 15587/108\ h^{8}s^{4}\\ 
\end{array}
\\
\begin{array}{rcl}
\sigma_{(3,5,6,3,4,3,1)} & = & 2401/54\ h^{24} - 7957/27\ h^{20}s - 245/3\ h^{18}t + 5249/9\ h^{16}s^{2} + 7562/27\ h^{14}st\\
 & &  - 10214/27\ h^{12}s^{3} + 14707/108\ h^{8}s^{4}\\ 
\end{array}
\\
\begin{array}{rcl}
\sigma_{(3,4,6,3,4,3,2)} & = & 1799/162\ h^{24} - 5363/81\ h^{20}s - 679/27\ h^{18}t + 2590/27\ h^{16}s^{2}\\
 & &  + 6934/81\ h^{14}st + 919/162\ h^{12}s^{3} - 7345/324\ h^{8}s^{4}\\ 
\end{array}
\vspace{.5cm}
$

\noindent
Schubert cells in degree $25$:

\noindent
$
\begin{array}{rcl}
\sigma_{(2,4,6,3,5,4,2)} & = & 7808/81\ h^{25} - 51742/81\ h^{21}s - 9685/54\ h^{19}t + 33914/27\ h^{17}s^{2}\\
 & &  + 100099/162\ h^{15}st - 128899/162\ h^{13}s^{3} + 91423/324\ h^{9}s^{4} + s^{2}q\\ 
\end{array}
\\
\begin{array}{rcl}
\sigma_{(3,5,6,3,5,3,1)} & = & -6019/108\ h^{25} + 37841/108\ h^{21}s + 465/4\ h^{19}t - 21511/36\ h^{17}s^{2}\\
 & &  - 10777/27\ h^{15}st + 10781/54\ h^{13}s^{3} - 1727/108\ h^{9}s^{4}\\ 
\end{array}
\\
\begin{array}{rcl}
\sigma_{(3,4,6,3,5,3,2)} & = & -14329/162\ h^{25} + 46990/81\ h^{21}s + 9073/54\ h^{19}t - 120239/108\ h^{17}s^{2}\\
 & &  - 93601/162\ h^{15}st + 105661/162\ h^{13}s^{3} - 17368/81\ h^{9}s^{4} - s^{2}q\\ 
\end{array}
\\
\begin{array}{rcl}
\sigma_{(3,5,6,3,4,3,2)} & = & 38279/324\ h^{25} - 247279/324\ h^{21}s - 25049/108\ h^{19}t + 76063/54\ h^{17}s^{2}\\
 & &  + 64517/81\ h^{15}st - 57767/81\ h^{13}s^{3} + 15964/81\ h^{9}s^{4} + s^{2}q\\ 
\end{array}
\vspace{.5cm}
$

\noindent
Schubert cells in degree $26$:

\noindent
$
\begin{array}{rcl}
\sigma_{(3,5,7,3,5,3,1)} & = & -23401/162\ h^{26} + 147137/162\ h^{22}s + 8132/27\ h^{20}t - 83665/54\ h^{18}s^{2}\\
 & &  - 167533/162\ h^{16}st + 41990/81\ h^{14}s^{3} - 6731/162\ h^{10}s^{4}\\ 
\end{array}
\\
\begin{array}{rcl}
\sigma_{(3,5,6,3,5,3,2)} & = & 929/162\ h^{26} - 16321/648\ h^{22}s - 2035/108\ h^{20}t - 2449/216\ h^{18}s^{2}\\
 & &  + 41231/648\ h^{16}st + 10196/81\ h^{14}s^{3} - 6589/81\ h^{10}s^{4}\\ 
\end{array}
\\
\begin{array}{rcl}
\sigma_{(3,4,6,3,5,4,2)} & = & 994/27\ h^{26} - 2995/12\ h^{22}s - 65\ h^{20}t + 2069/4\ h^{18}s^{2} + 2695/12\ h^{16}st\\
 & &  - 20483/54\ h^{14}s^{3} + 5425/36\ h^{10}s^{4}\\ 
\end{array}
\vspace{.5cm}
$

\noindent
Schubert cells in degree $27$:

\noindent
$
\begin{array}{rcl}
\sigma_{(3,5,7,4,5,3,1)} & = & -22028/81\ h^{27} + 138517/81\ h^{23}s + 15305/27\ h^{21}t - 78782/27\ h^{19}s^{2}\\
 & &  - 315313/162\ h^{17}st + 79169/81\ h^{15}s^{3} - 6349/81\ h^{11}s^{4} + 4\ stq\\ 
\end{array}
\\
\begin{array}{rcl}
\sigma_{(3,5,7,3,5,3,2)} & = & -205/2\ h^{27} + 8281/12\ h^{23}s + 551/3\ h^{21}t - 1409\ h^{19}s^{2} - 22811/36\ h^{17}st\\
 & &  + 1989/2\ h^{15}s^{3} - 1153/3\ h^{11}s^{4} - 2\ stq\\ 
\end{array}
\\
\begin{array}{rcl}
\sigma_{(3,5,6,3,5,4,2)} & = & 8767/81\ h^{27} - 463495/648\ h^{23}s - 21871/108\ h^{21}t + 301895/216\ h^{19}s^{2}\\
 & &  + 451829/648\ h^{17}st - 140717/162\ h^{15}s^{3} + 24542/81\ h^{11}s^{4} + stq\\ 
\end{array}
\vspace{.5cm}
$

\noindent
Schubert cells in degree $28$:

\noindent
$
\begin{array}{rcl}
\sigma_{(3,5,7,4,5,3,2)} & = & -218447/486\ h^{28} + 710360/243\ h^{24}s + 141325/162\ h^{22}t - 888215/162\ h^{20}s^{2}\\
 & &  - 1458431/486\ h^{18}st + 2871745/972\ h^{16}s^{3} - 431527/486\ h^{12}s^{4}\\ 
\end{array}
\\
\begin{array}{rcl}
\sigma_{(3,5,7,3,5,4,2)} & = & 107965/486\ h^{28} - 703829/486\ h^{24}s - 34714/81\ h^{22}t + 442939/162\ h^{20}s^{2}\\
 & &  + 716365/486\ h^{18}st - 367151/243\ h^{16}s^{3} + 454465/972\ h^{12}s^{4}\\ 
\end{array}
\vspace{.5cm}
$

\noindent
Schubert cells in degree $29$:

\noindent
$
\begin{array}{rcl}
\sigma_{(3,5,7,4,5,4,2)} & = & -78557/108\ h^{29} + 14254/3\ h^{25}s + 452857/324\ h^{23}t - 36057/4\ h^{21}s^{2}\\
 & &  - 519361/108\ h^{19}st + 136862/27\ h^{17}s^{3} - 14396/9\ h^{13}s^{4} - 2/3\ s^{3}q\\ 
\end{array}
\\
\begin{array}{rcl}
\sigma_{(3,5,7,3,6,4,2)} & = & 179729/162\ h^{29} - 4705735/648\ h^{25}s - 229097/108\ h^{23}t + 2992283/216\ h^{21}s^{2}\\
 & &  + 4729385/648\ h^{19}st - 2574205/324\ h^{17}s^{3} + 829661/324\ h^{13}s^{4} + s^{3}q\\ 
\end{array}
\vspace{.5cm}
$

\noindent
Schubert cells in degree $30$:

\noindent
$
\begin{array}{rcl}
\sigma_{(3,5,7,4,6,4,2)} & = & 165809/1944\ h^{30} - 140594/243\ h^{26}s - 12083/81\ h^{24}t + 97579/81\ h^{22}s^{2}\\
 & &  + 124891/243\ h^{20}st - 217618/243\ h^{18}s^{3} + 87011/243\ h^{14}s^{4}\\ 
\end{array}
\vspace{.5cm}
$

\noindent
Schubert cells in degree $31$:

\noindent
$
\begin{array}{rcl}
\sigma_{(3,5,8,4,6,4,2)} & = & 353291/729\ h^{31} - 4707605/1458\ h^{27}s - 215872/243\ h^{25}t\\
 & &  + 6258959/972\ h^{23}s^{2} + 4460419/1458\ h^{21}st - 6209929/1458\ h^{19}s^{3}\\
 & &  + 2281865/1458\ h^{15}s^{4}\\ 
\end{array}
\vspace{.5cm}
$

\noindent
Schubert cells in degree $32$:

\noindent
$
\begin{array}{rcl}
\sigma_{(3,6,8,4,6,4,2)} & = & 1302991/729\ h^{32} - 8617946/729\ h^{28}s - 810560/243\ h^{26}t\\
 & &  + 5626099/243\ h^{24}s^{2} + 8371429/729\ h^{22}st - 10566061/729\ h^{20}s^{3}\\
 & &  + 29652913/5832\ h^{16}s^{4}\\ 
\end{array}
\vspace{.5cm}
$

\noindent
Schubert cells in degree $33$:

\noindent
$
\begin{array}{rcl}
\sigma_{(4,6,8,4,6,4,2)} & = & 15615181/2916\ h^{33} - 205589299/5832\ h^{29}s - 9822617/972\ h^{27}t\\
 & &  + 132640073/1944\ h^{25}s^{2} + 202845737/5832\ h^{23}st - 119894383/2916\ h^{21}s^{3}\\
 & &  + 40610093/2916\ h^{17}s^{4} + 1/3\ s^{4}q\\ 
\end{array}
\vspace{.5cm}
$

\subsection{$E_7/P_7$}


\noindent
Schubert cells in degree $0$:

\noindent
$
\begin{array}{rcl}
\sigma_{(0,0,0,0,0,0,0)} & = & 1\ \\ 
\end{array}
\vspace{.5cm}
$

\noindent
Schubert cells in degree $1$:

\noindent
$
\begin{array}{rcl}
\sigma_{(0,0,0,0,0,0,1)} & = & h\\ 
\end{array}
\vspace{.5cm}
$

\noindent
Schubert cells in degree $2$:

\noindent
$
\begin{array}{rcl}
\sigma_{(0,0,0,0,0,1,1)} & = & h^{2}\\ 
\end{array}
\vspace{.5cm}
$

\noindent
Schubert cells in degree $3$:

\noindent
$
\begin{array}{rcl}
\sigma_{(0,0,0,0,1,1,1)} & = & h^{3}\\ 
\end{array}
\vspace{.5cm}
$

\noindent
Schubert cells in degree $4$:

\noindent
$
\begin{array}{rcl}
\sigma_{(0,0,0,1,1,1,1)} & = & h^{4}\\ 
\end{array}
\vspace{.5cm}
$

\noindent
Schubert cells in degree $5$:

\noindent
$
\begin{array}{rcl}
\sigma_{(0,0,1,1,1,1,1)} & = & h^{5} - s\\ 
\end{array}
\\
\begin{array}{rcl}
\sigma_{(0,1,0,1,1,1,1)} & = & s\\ 
\end{array}
\vspace{.5cm}
$

\noindent
Schubert cells in degree $6$:

\noindent
$
\begin{array}{rcl}
\sigma_{(1,0,1,1,1,1,1)} & = & h^{6} - 2\ hs\\ 
\end{array}
\\
\begin{array}{rcl}
\sigma_{(0,1,1,1,1,1,1)} & = & hs\\ 
\end{array}
\vspace{.5cm}
$

\noindent
Schubert cells in degree $7$:

\noindent
$
\begin{array}{rcl}
\sigma_{(1,1,1,1,1,1,1)} & = & h^{7} - 2\ h^{2}s\\ 
\end{array}
\\
\begin{array}{rcl}
\sigma_{(0,1,1,2,1,1,1)} & = & -h^{7} + 3\ h^{2}s\\ 
\end{array}
\vspace{.5cm}
$

\noindent
Schubert cells in degree $8$:

\noindent
$
\begin{array}{rcl}
\sigma_{(1,1,1,2,1,1,1)} & = & h^{8} - 2\ h^{3}s\\ 
\end{array}
\\
\begin{array}{rcl}
\sigma_{(0,1,1,2,2,1,1)} & = & -2\ h^{8} + 5\ h^{3}s\\ 
\end{array}
\vspace{.5cm}
$

\noindent
Schubert cells in degree $9$:

\noindent
$
\begin{array}{rcl}
\sigma_{(1,1,2,2,1,1,1)} & = & t\\ 
\end{array}
\\
\begin{array}{rcl}
\sigma_{(1,1,1,2,2,1,1)} & = & h^{9} - 2\ h^{4}s - t\\ 
\end{array}
\\
\begin{array}{rcl}
\sigma_{(0,1,1,2,2,2,1)} & = & -3\ h^{9} + 7\ h^{4}s + t\\ 
\end{array}
\vspace{.5cm}
$

\noindent
Schubert cells in degree $10$:

\noindent
$
\begin{array}{rcl}
\sigma_{(1,1,2,2,2,1,1)} & = & ht\\ 
\end{array}
\\
\begin{array}{rcl}
\sigma_{(1,1,1,2,2,2,1)} & = & h^{10} - 2\ h^{5}s - 2\ ht\\ 
\end{array}
\\
\begin{array}{rcl}
\sigma_{(0,1,1,2,2,2,2)} & = & -4\ h^{10} + 9\ h^{5}s + 3\ ht\\ 
\end{array}
\vspace{.5cm}
$

\noindent
Schubert cells in degree $11$:

\noindent
$
\begin{array}{rcl}
\sigma_{(1,1,2,3,2,1,1)} & = & -5\ h^{11} + 11\ h^{6}s + 6\ h^{2}t\\ 
\end{array}
\\
\begin{array}{rcl}
\sigma_{(1,1,2,2,2,2,1)} & = & 5\ h^{11} - 11\ h^{6}s - 5\ h^{2}t\\ 
\end{array}
\\
\begin{array}{rcl}
\sigma_{(1,1,1,2,2,2,2)} & = & -4\ h^{11} + 9\ h^{6}s + 3\ h^{2}t\\ 
\end{array}
\vspace{.5cm}
$

\noindent
Schubert cells in degree $12$:

\noindent
$
\begin{array}{rcl}
\sigma_{(1,2,2,3,2,1,1)} & = & -14\ h^{12} + 31\ h^{7}s + 14\ h^{3}t\\ 
\end{array}
\\
\begin{array}{rcl}
\sigma_{(1,1,2,3,2,2,1)} & = & 9\ h^{12} - 20\ h^{7}s - 8\ h^{3}t\\ 
\end{array}
\\
\begin{array}{rcl}
\sigma_{(1,1,2,2,2,2,2)} & = & -4\ h^{12} + 9\ h^{7}s + 3\ h^{3}t\\ 
\end{array}
\vspace{.5cm}
$

\noindent
Schubert cells in degree $13$:

\noindent
$
\begin{array}{rcl}
\sigma_{(1,2,2,3,2,2,1)} & = & -14\ h^{13} + 31\ h^{8}s + 14\ h^{4}t\\ 
\end{array}
\\
\begin{array}{rcl}
\sigma_{(1,1,2,3,3,2,1)} & = & 27\ h^{13} - 60\ h^{8}s - 25\ h^{4}t\\ 
\end{array}
\\
\begin{array}{rcl}
\sigma_{(1,1,2,3,2,2,2)} & = & -4\ h^{13} + 9\ h^{8}s + 3\ h^{4}t\\ 
\end{array}
\vspace{.5cm}
$

\noindent
Schubert cells in degree $14$:

\noindent
$
\begin{array}{rcl}
\sigma_{(1,2,2,3,3,2,1)} & = & 17/2\ h^{14} - 19\ h^{9}s - 7\ h^{5}t\\ 
\end{array}
\\
\begin{array}{rcl}
\sigma_{(1,2,2,3,2,2,2)} & = & -45/2\ h^{14} + 50\ h^{9}s + 21\ h^{5}t\\ 
\end{array}
\\
\begin{array}{rcl}
\sigma_{(1,1,2,3,3,2,2)} & = & 37/2\ h^{14} - 41\ h^{9}s - 18\ h^{5}t\\ 
\end{array}
\vspace{.5cm}
$

\noindent
Schubert cells in degree $15$:

\noindent
$
\begin{array}{rcl}
\sigma_{(1,2,2,4,3,2,1)} & = & 31\ h^{15} - 69\ h^{10}s - 28\ h^{6}t\\ 
\end{array}
\\
\begin{array}{rcl}
\sigma_{(1,2,2,3,3,2,2)} & = & -45/2\ h^{15} + 50\ h^{10}s + 21\ h^{6}t\\ 
\end{array}
\\
\begin{array}{rcl}
\sigma_{(1,1,2,3,3,3,2)} & = & 41\ h^{15} - 91\ h^{10}s - 39\ h^{6}t\\ 
\end{array}
\vspace{.5cm}
$

\noindent
Schubert cells in degree $16$:

\noindent
$
\begin{array}{rcl}
\sigma_{(1,2,3,4,3,2,1)} & = & 189/2\ h^{16} - 210\ h^{11}s - 88\ h^{7}t\\ 
\end{array}
\\
\begin{array}{rcl}
\sigma_{(1,2,2,4,3,2,2)} & = & -127/2\ h^{16} + 141\ h^{11}s + 60\ h^{7}t\\ 
\end{array}
\\
\begin{array}{rcl}
\sigma_{(1,2,2,3,3,3,2)} & = & 41\ h^{16} - 91\ h^{11}s - 39\ h^{7}t\\ 
\end{array}
\vspace{.5cm}
$

\noindent
Schubert cells in degree $17$:

\noindent
$
\begin{array}{rcl}
\sigma_{(2,2,3,4,3,2,1)} & = & 199\ h^{17} - 442\ h^{12}s - 187\ h^{8}t\\ 
\end{array}
\\
\begin{array}{rcl}
\sigma_{(1,2,3,4,3,2,2)} & = & -209/2\ h^{17} + 232\ h^{12}s + 99\ h^{8}t\\ 
\end{array}
\\
\begin{array}{rcl}
\sigma_{(1,2,2,4,3,3,2)} & = & 41\ h^{17} - 91\ h^{12}s - 39\ h^{8}t\\ 
\end{array}
\vspace{.5cm}
$

\noindent
Schubert cells in degree $18$:

\noindent
$
\begin{array}{rcl}
\sigma_{(2,2,3,4,3,2,2)} & = & 72\ h^{18} - 160\ h^{13}s - 67\ h^{9}t - t^{2}\\ 
\end{array}
\\
\begin{array}{rcl}
\sigma_{(1,2,3,4,3,3,2)} & = & -353/2\ h^{18} + 392\ h^{13}s + 166\ h^{9}t + t^{2}\\ 
\end{array}
\\
\begin{array}{rcl}
\sigma_{(1,2,2,4,4,3,2)} & = & 435/2\ h^{18} - 483\ h^{13}s - 205\ h^{9}t - t^{2}\\ 
\end{array}
\vspace{.5cm}
$

\noindent
Schubert cells in degree $19$:

\noindent
$
\begin{array}{rcl}
\sigma_{(2,2,3,4,3,3,2)} & = & -449/2\ h^{19} + 997/2\ h^{14}s + 212\ h^{10}t\\ 
\end{array}
\\
\begin{array}{rcl}
\sigma_{(1,2,3,4,4,3,2)} & = & 531/4\ h^{19} - 1179/4\ h^{14}s - 251/2\ h^{10}t\\ 
\end{array}
\vspace{.5cm}
$

\noindent
Schubert cells in degree $20$:

\noindent
$
\begin{array}{rcl}
\sigma_{(2,2,3,4,4,3,2)} & = & -1745/4\ h^{20} + 3875/4\ h^{15}s + 823/2\ h^{11}t\\ 
\end{array}
\\
\begin{array}{rcl}
\sigma_{(1,2,3,5,4,3,2)} & = & 569\ h^{20} - 2527/2\ h^{15}s - 537\ h^{11}t\\ 
\end{array}
\vspace{.5cm}
$

\noindent
Schubert cells in degree $21$:

\noindent
$
\begin{array}{rcl}
\sigma_{(2,2,3,5,4,3,2)} & = & -648\ h^{21} + 1439\ h^{16}s + 611\ h^{12}t\\ 
\end{array}
\\
\begin{array}{rcl}
\sigma_{(1,3,3,5,4,3,2)} & = & 1217\ h^{21} - 5405/2\ h^{16}s - 1148\ h^{12}t\\ 
\end{array}
\vspace{.5cm}
$

\noindent
Schubert cells in degree $22$:

\noindent
$
\begin{array}{rcl}
\sigma_{(2,2,4,5,4,3,2)} & = & -8307/4\ h^{22} + 18447/4\ h^{17}s + 3917/2\ h^{13}t\\ 
\end{array}
\\
\begin{array}{rcl}
\sigma_{(2,3,3,5,4,3,2)} & = & 1217\ h^{22} - 5405/2\ h^{17}s - 1148\ h^{13}t\\ 
\end{array}
\vspace{.5cm}
$

\noindent
Schubert cells in degree $23$:

\noindent
$
\begin{array}{rcl}
\sigma_{(2,3,4,5,4,3,2)} & = & -2143/4\ h^{23} + 4759/4\ h^{18}s + 505\ h^{14}t\\ 
\end{array}
\vspace{.5cm}
$

\noindent
Schubert cells in degree $24$:

\noindent
$
\begin{array}{rcl}
\sigma_{(2,3,4,6,4,3,2)} & = & -4577/2\ h^{24} + 5082\ h^{19}s + 2158\ h^{15}t\\ 
\end{array}
\vspace{.5cm}
$

\noindent
Schubert cells in degree $25$:

\noindent
$
\begin{array}{rcl}
\sigma_{(2,3,4,6,5,3,2)} & = & -7335\ h^{25} + 32577/2\ h^{20}s + 13835/2\ h^{16}t\\ 
\end{array}
\vspace{.5cm}
$

\noindent
Schubert cells in degree $26$:

\noindent
$
\begin{array}{rcl}
\sigma_{(2,3,4,6,5,4,2)} & = & -62701/4\ h^{26} + 139237/4\ h^{21}s + 29567/2\ h^{17}t\\ 
\end{array}
\vspace{.5cm}
$

\noindent
Schubert cells in degree $27$:

\noindent
$
\begin{array}{rcl}
\sigma_{(2,3,4,6,5,4,3)} & = & -83785/4\ h^{27} + 186057/4\ h^{22}s + 19755\ h^{18}t - t^{3}\\ 
\end{array}
\vspace{.5cm}
$

\newpage

\section{$E_8/P_8$}

\subsection{Giambelli formulas}

\input{e8p8}

\newpage

\subsection{Tableaux}

Here are the 9 tableaux involved to compute $s^2$: 
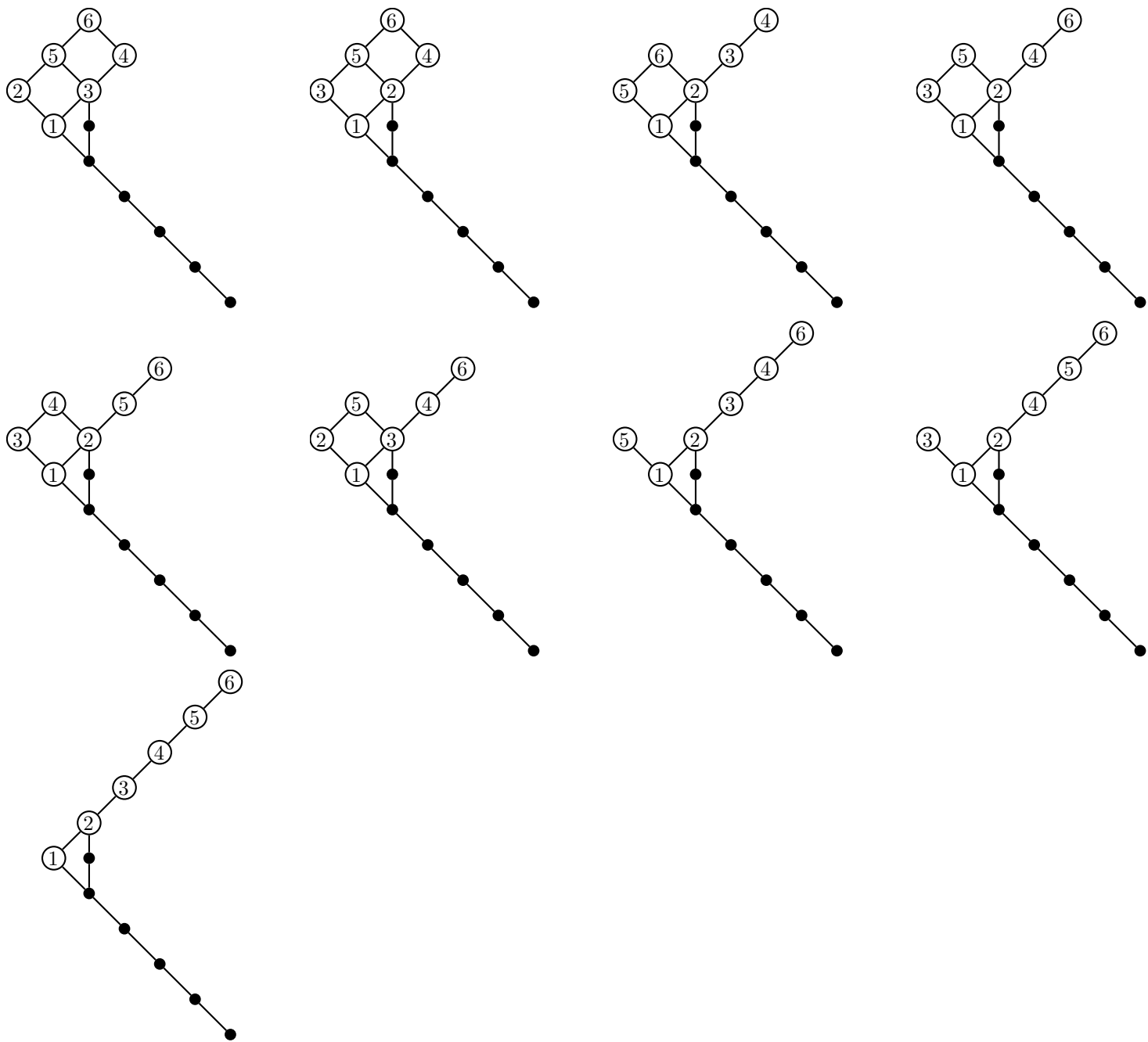
\vspace{-.5cm}

$$
\hspace{-1.5cm}
\begin{array}{cccccc}
\input{tableau-e8p8-12-1} & \input{tableau-e8p8-12-2} & \input{tableau-e8p8-12-3} & \input{tableau-e8p8-12-4} \\
\input{tableau-e8p8-12-5} & \input{tableau-e8p8-12-6} & \input{tableau-e8p8-12-7} & \input{tableau-e8p8-12-8} \\
\input{tableau-e8p8-12-9}
\end{array}
$$

\newpage

Here are the 21 tableaux involved to compute $s\cdot t$:

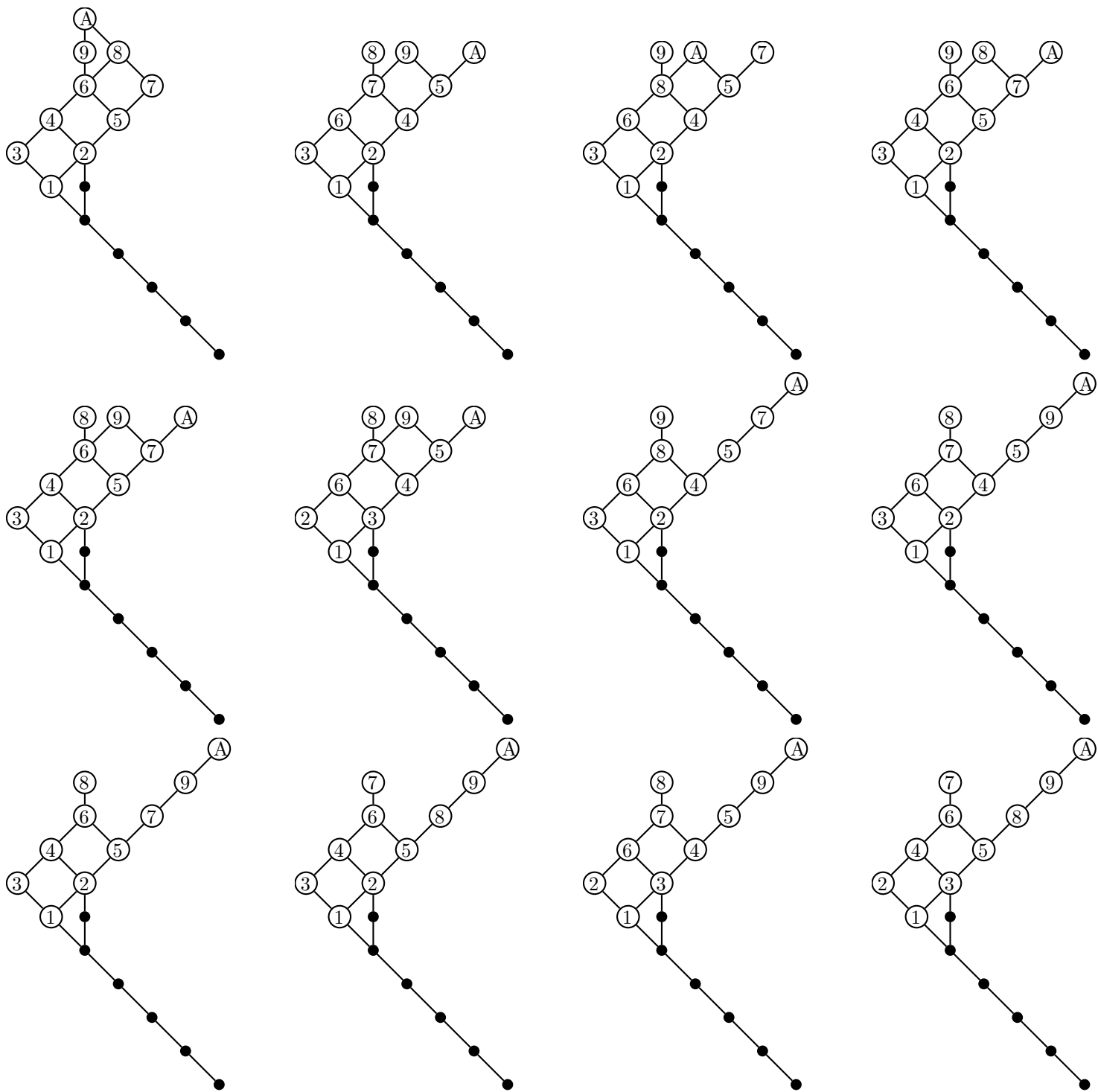
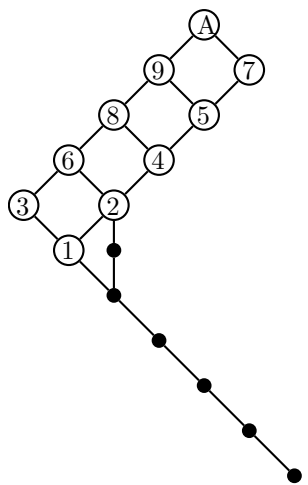
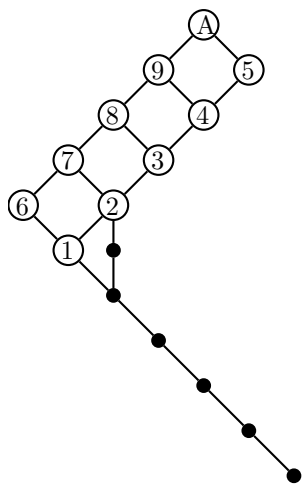
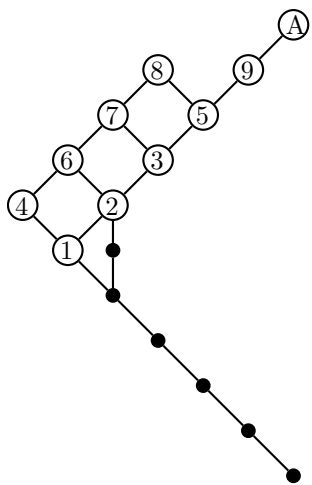
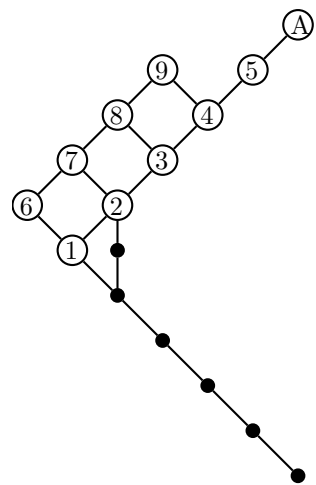
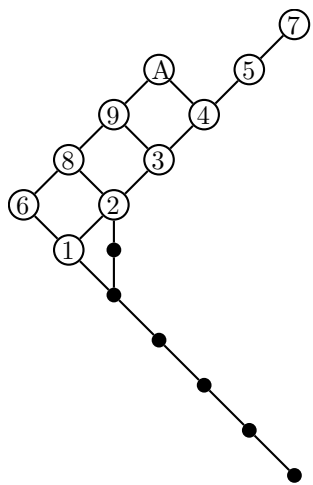
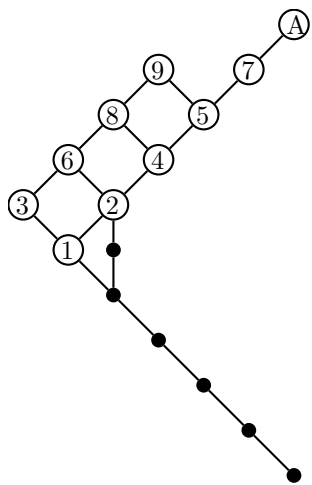
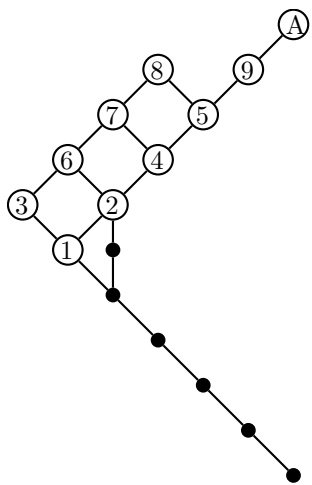
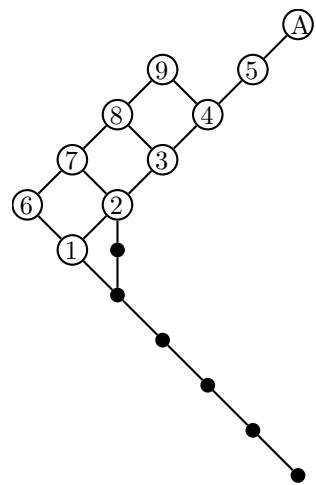
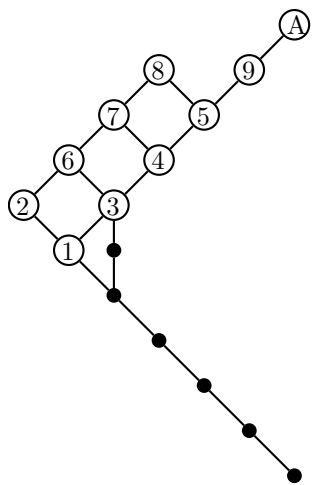
$$
\hspace{-1.5cm}
\begin{array}{cccccc}
\input{tableau-e8p8-16-1} & \input{tableau-e8p8-16-2} & \input{tableau-e8p8-16-3} & \input{tableau-e8p8-16-4} \\
\input{tableau-e8p8-16-5} & \input{tableau-e8p8-16-6} & \input{tableau-e8p8-16-7} & \input{tableau-e8p8-16-8} \\
\input{tableau-e8p8-16-9} & \input{tableau-e8p8-16-10} & \input{tableau-e8p8-16-11} & \input{tableau-e8p8-16-12}
\end{array}
$$

$$
\hspace{-1.5cm}
\begin{array}{cccccc}
\input{tableau-e8p8-16-13} & \input{tableau-e8p8-16-14} & \input{tableau-e8p8-16-15} & \input{tableau-e8p8-16-16} \\
\input{tableau-e8p8-16-17} & \input{tableau-e8p8-16-18} & \input{tableau-e8p8-16-19} & \input{tableau-e8p8-16-16} \\
\input{tableau-e8p8-16-21}
\end{array}
$$

\newpage

Here are the 3 tableaux involved to compute $s \cdot \s \poidsehuit 01222221$:

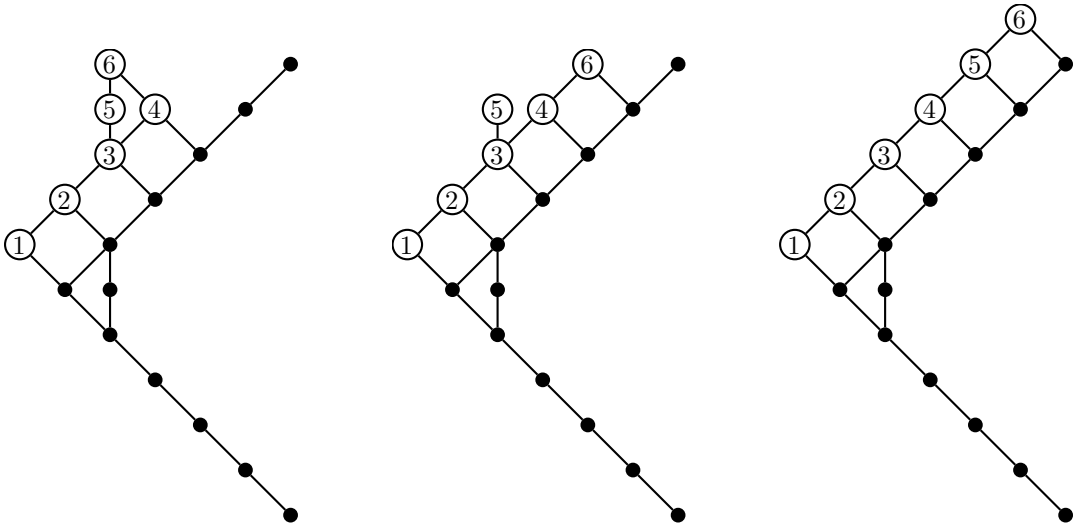
$$
\begin{array}{cccccc}
\input{tableau-e8p8-18-1} & \input{tableau-e8p8-18-2} & \input{tableau-e8p8-18-3}
\end{array}
$$

\newpage

Here are the 51 tableaux involved to compute $t^2$: 
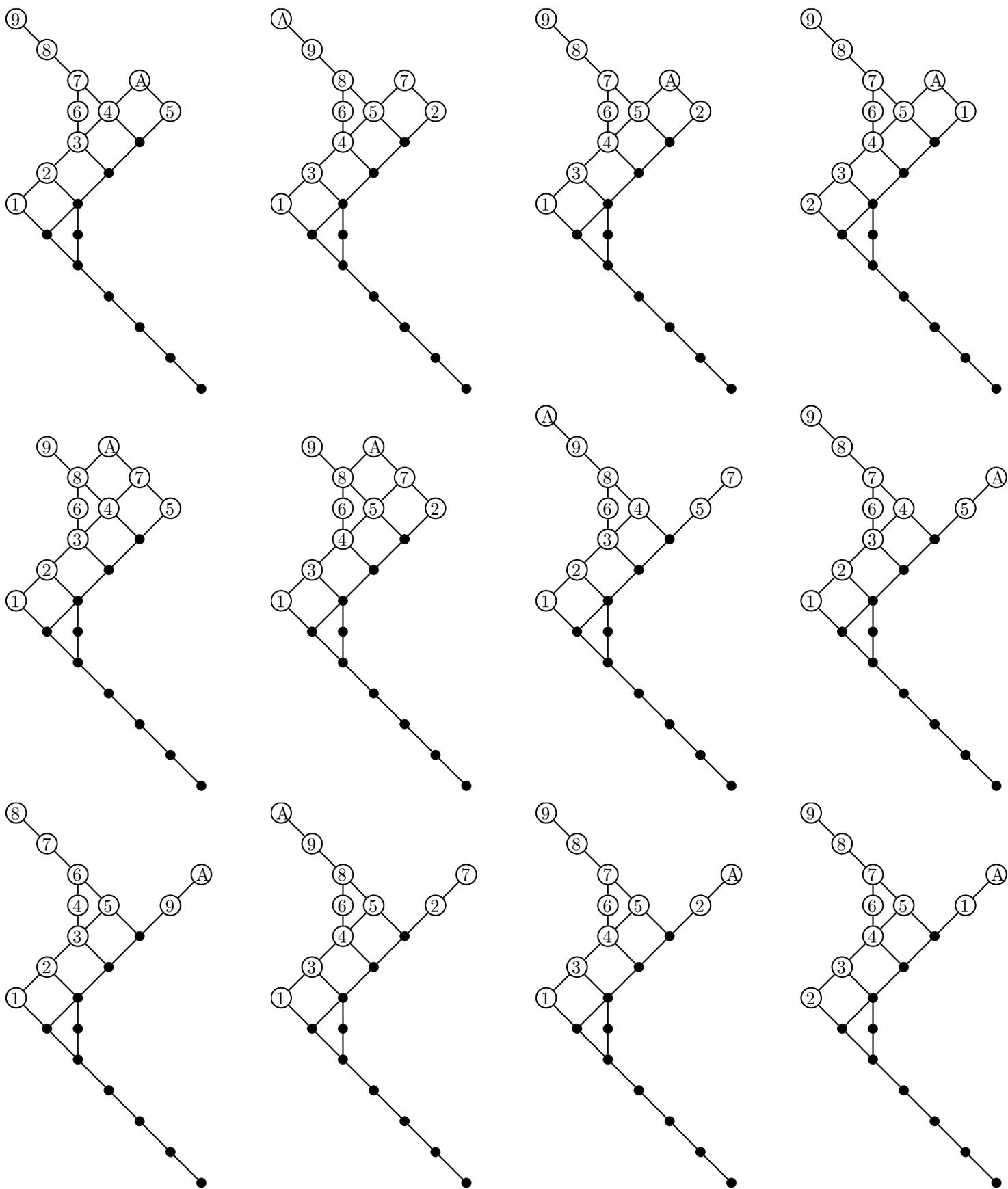
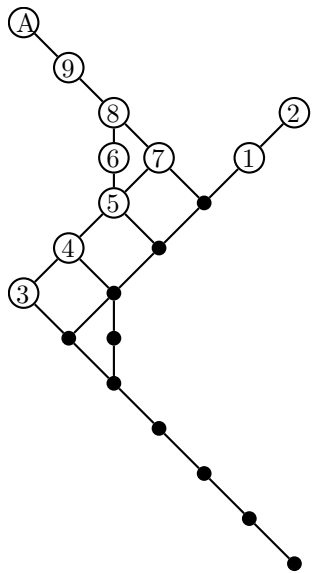
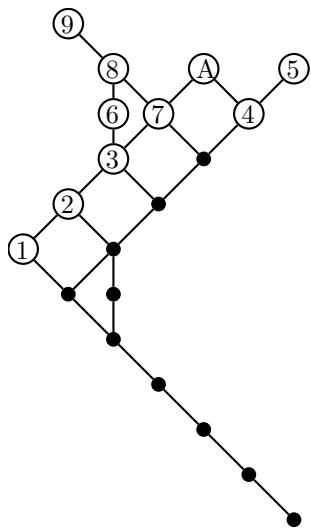
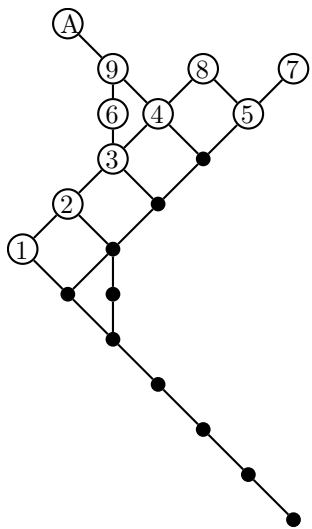
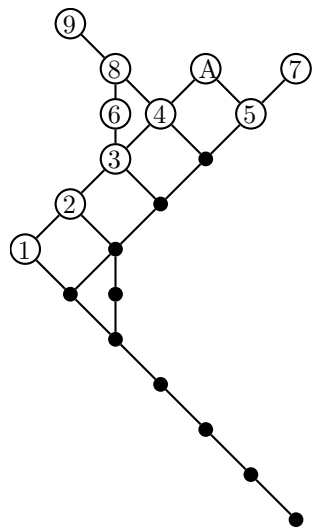
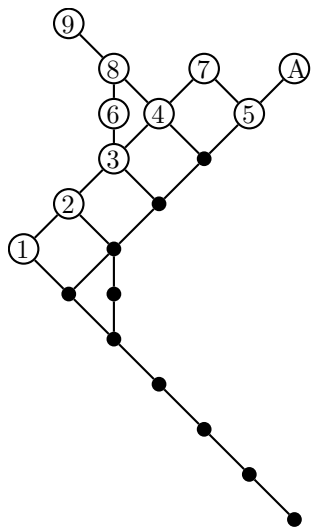
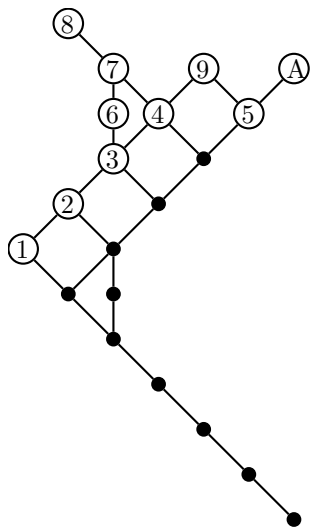
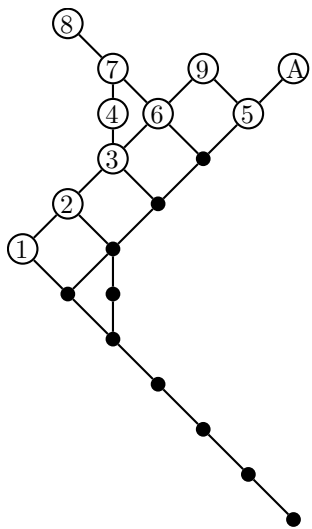
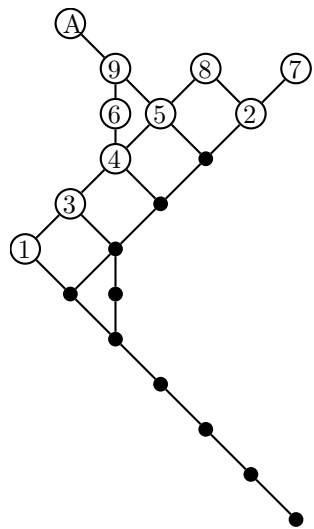
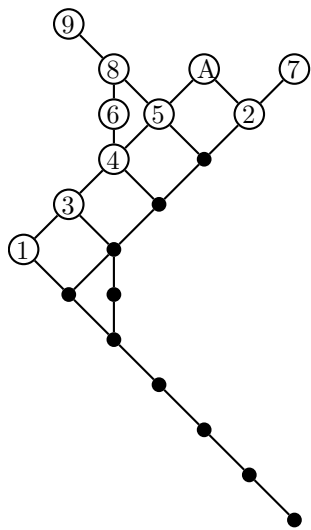
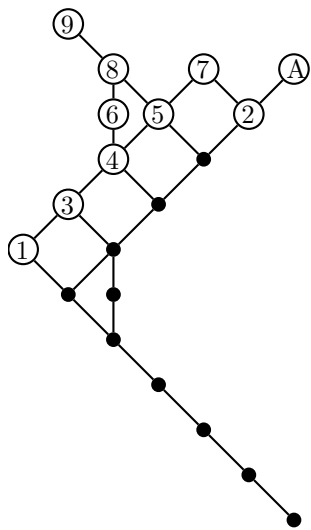
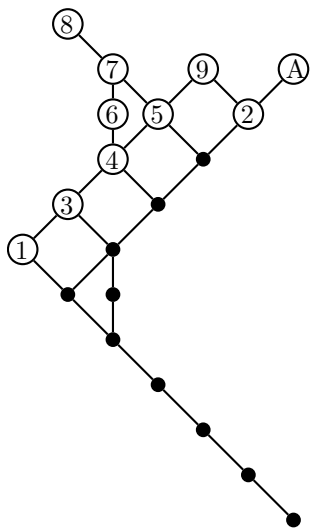
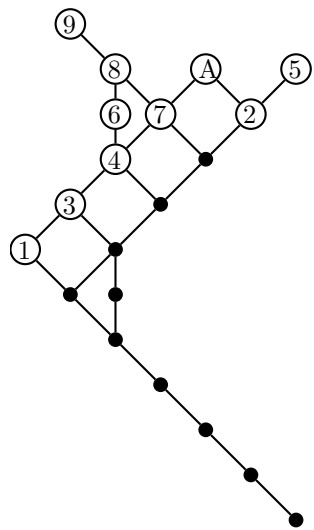
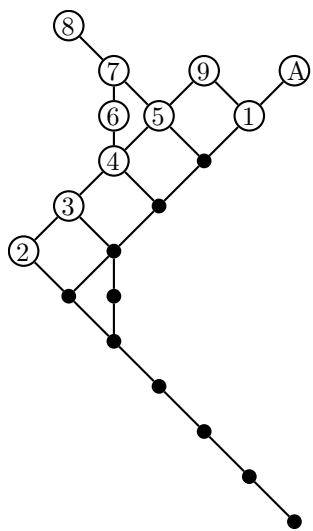
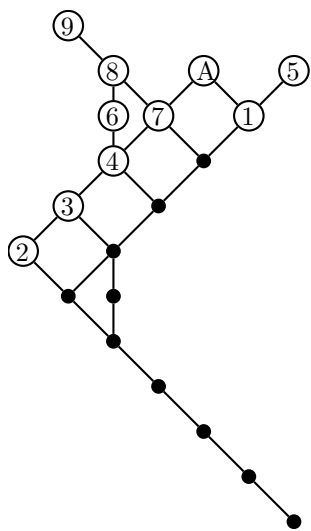
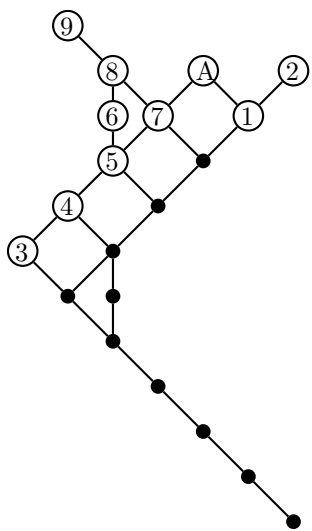
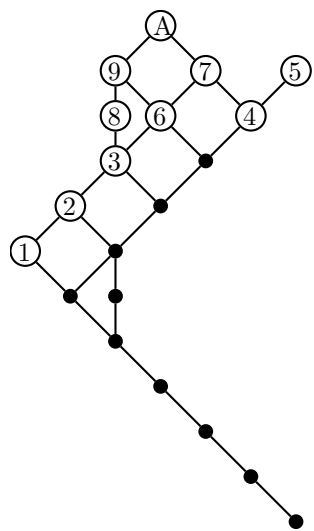
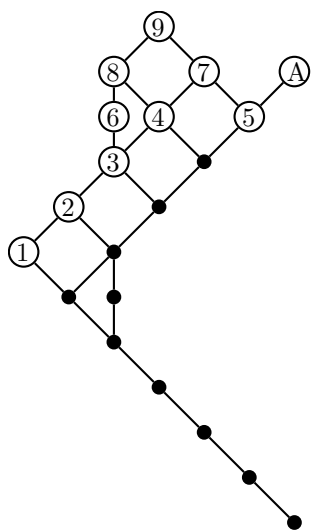
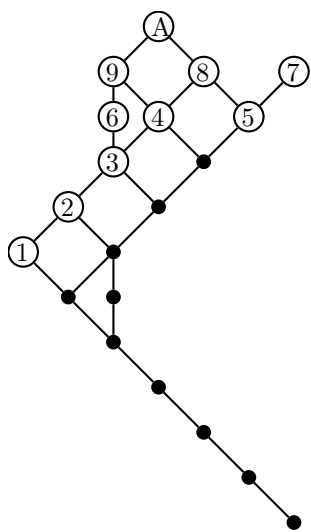
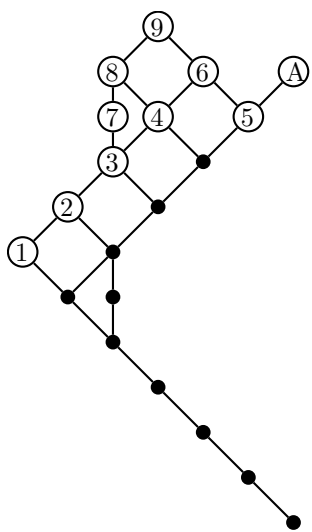
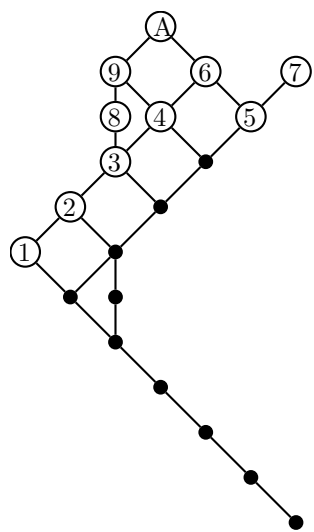
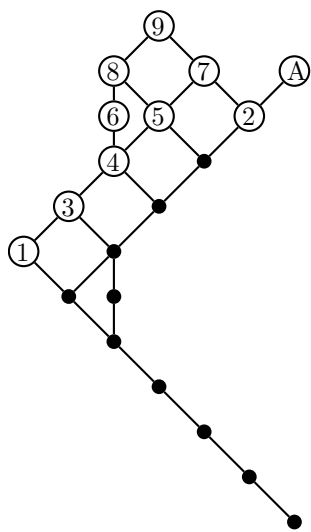
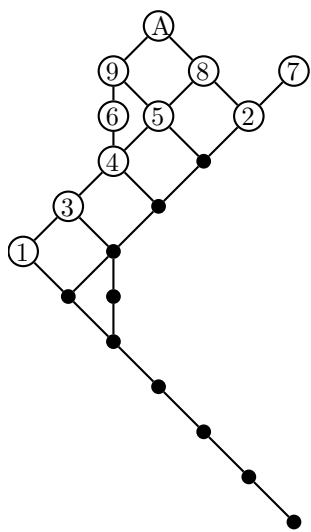
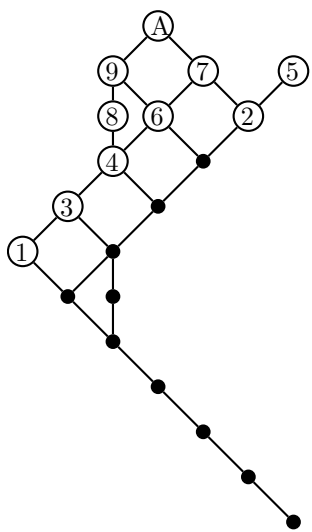
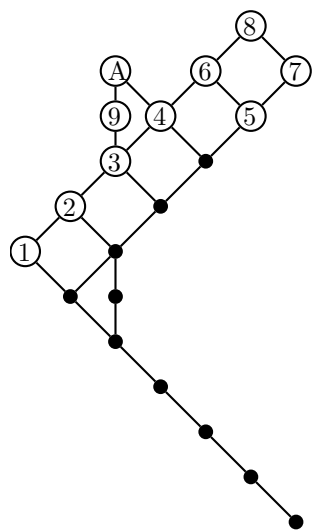
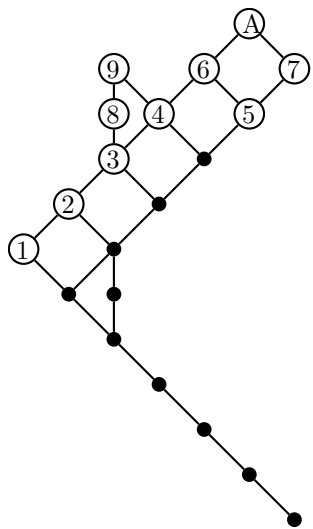
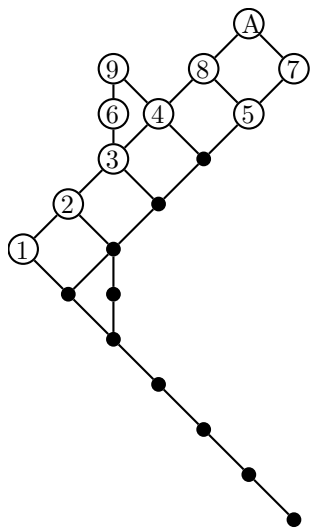
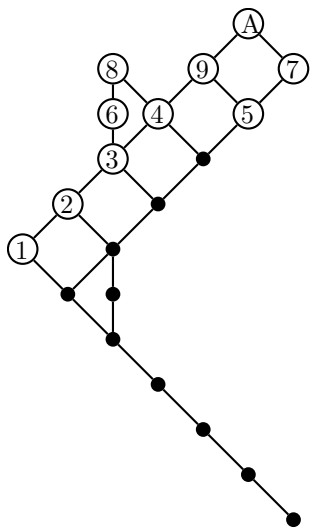
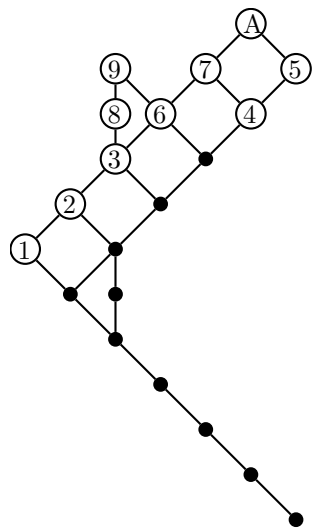
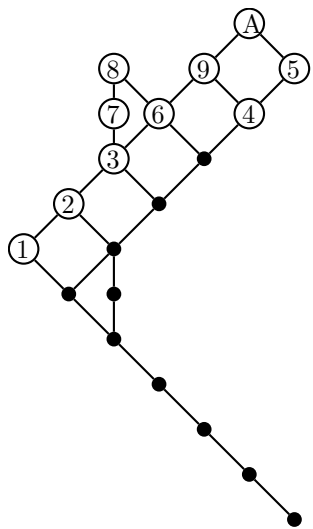
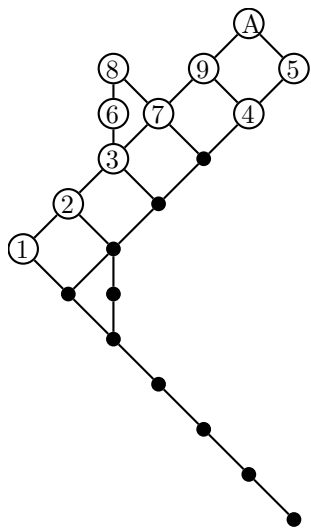
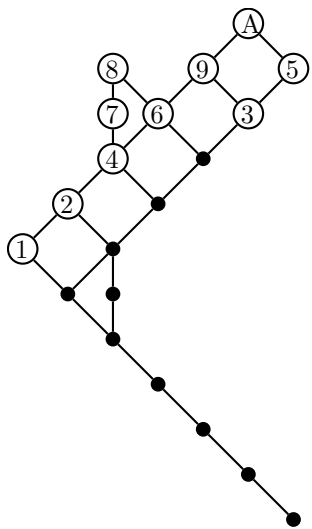
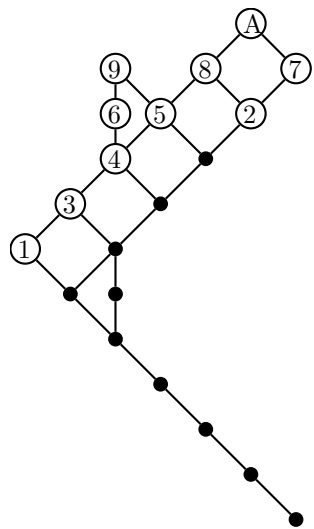
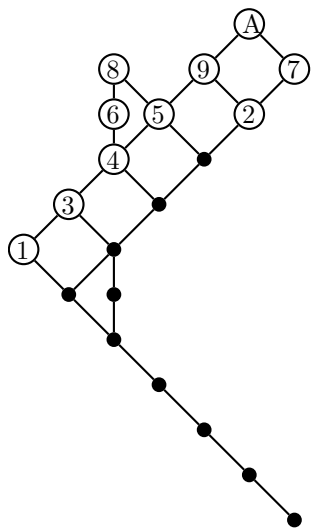
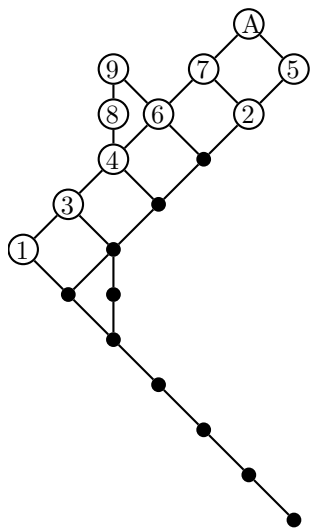
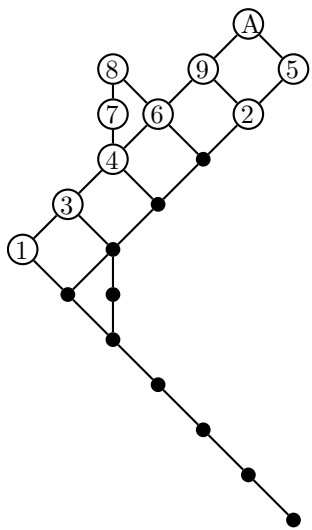
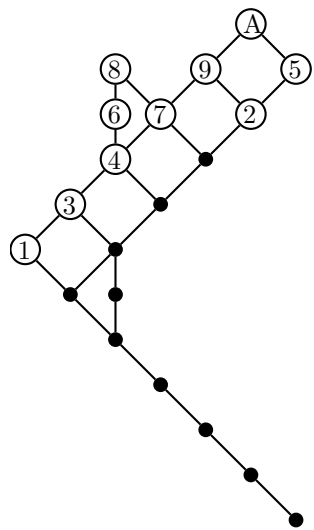
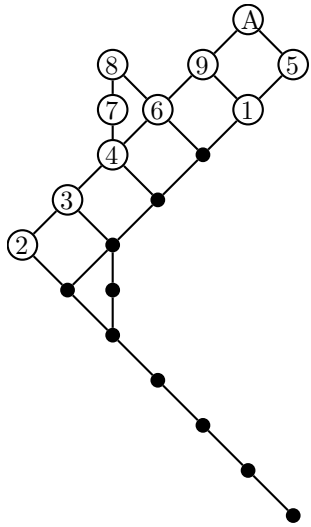
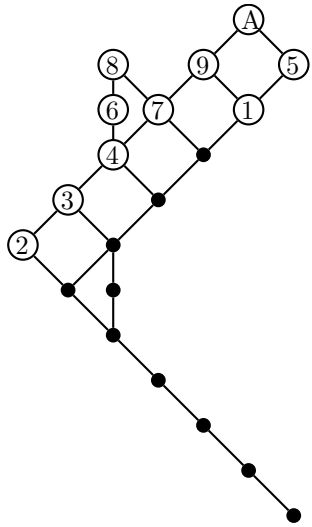
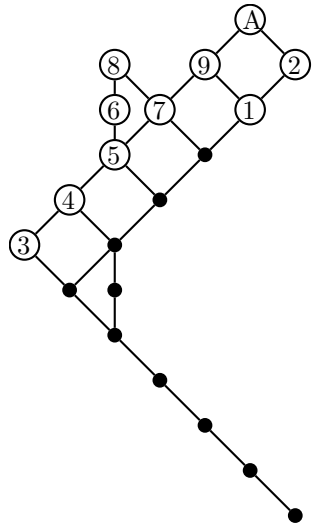
\vspace{-.5cm}

$$
\hspace{-1.5cm}
\begin{array}{cccccc}
\input{tableau-e8p8-20-1} & \input{tableau-e8p8-20-2} & \input{tableau-e8p8-20-3} & \input{tableau-e8p8-20-4} \\
\input{tableau-e8p8-20-5} & \input{tableau-e8p8-20-6} & \input{tableau-e8p8-20-7} & \input{tableau-e8p8-20-8} \\
\input{tableau-e8p8-20-9} & \input{tableau-e8p8-20-10} & \input{tableau-e8p8-20-11} & \input{tableau-e8p8-20-12}
\end{array}
$$

$$
\hspace{-1.5cm}
\begin{array}{cccccc}
\input{tableau-e8p8-20-13} & \input{tableau-e8p8-20-14} & \input{tableau-e8p8-20-15} & \input{tableau-e8p8-20-16} \\
\input{tableau-e8p8-20-17} & \input{tableau-e8p8-20-18} & \input{tableau-e8p8-20-19} & \input{tableau-e8p8-20-20} \\
\input{tableau-e8p8-20-21} & \input{tableau-e8p8-20-22} & \input{tableau-e8p8-20-23} & \input{tableau-e8p8-20-24}
\end{array}
$$

$$
\hspace{-1.5cm}
\begin{array}{cccccc}
\input{tableau-e8p8-20-25} & \input{tableau-e8p8-20-26} & \input{tableau-e8p8-20-27} & \input{tableau-e8p8-20-28} \\
\input{tableau-e8p8-20-29} & \input{tableau-e8p8-20-30} & \input{tableau-e8p8-20-31} & \input{tableau-e8p8-20-32} \\
\input{tableau-e8p8-20-33} & \input{tableau-e8p8-20-34} & \input{tableau-e8p8-20-35} & \input{tableau-e8p8-20-36}
\end{array}
$$

$$
\hspace{-1.5cm}
\begin{array}{cccccc}
\input{tableau-e8p8-20-37} & \input{tableau-e8p8-20-38} & \input{tableau-e8p8-20-39} & \input{tableau-e8p8-20-40} \\
\input{tableau-e8p8-20-41} & \input{tableau-e8p8-20-42} & \input{tableau-e8p8-20-43} & \input{tableau-e8p8-20-44} \\
\input{tableau-e8p8-20-45} & \input{tableau-e8p8-20-46} & \input{tableau-e8p8-20-47} & \input{tableau-e8p8-20-48}
\end{array}
$$

$$
\hspace{-1.5cm}
\begin{array}{cccccc} 
\input{tableau-e8p8-20-49} & \input{tableau-e8p8-20-50} & \input{tableau-e8p8-20-51} \\
\end{array}
$$

All of them rectify on

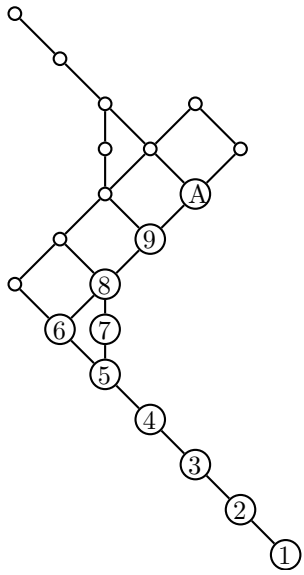
\centerline{ \input{reference-e8p8-20-0} }

\newpage

Here are the 7 tableaux involved to compute $t \cdot \s \poidsehuit 01222221$:

$$
\hspace{-1.5cm}
\begin{array}{cccccc}
\input{tableau-e8p8-22-1} & \input{tableau-e8p8-22-2} & \input{tableau-e8p8-22-3} & \input{tableau-e8p8-22-4} \\
\input{tableau-e8p8-22-5} & \input{tableau-e8p8-22-6} & \input{tableau-e8p8-22-7}
\end{array}
$$

\newpage

We have $\s \poidsehuit 01222221 ^2 = \s \poidsehuit 23444322$, in view of the following tableau:

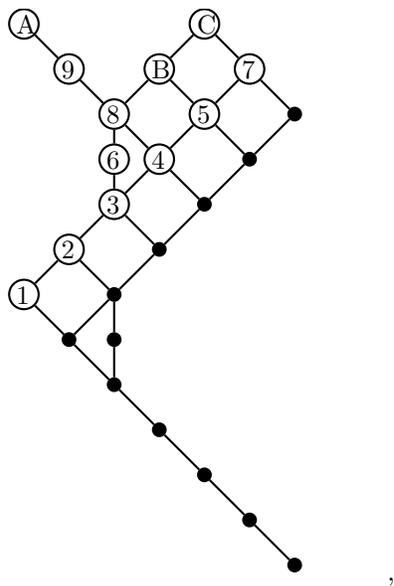
\centerline{
\input{tableau-e8p8-24-1}\ ,
}

\noindent
and the product $\s \poidsehuit 01222221 \cdot \s \poidsehuit 12332221$ is computed thanks to the 7 following tableaux:

$$
\hspace{-1.5cm}
\begin{array}{cccccc}
\input{tableau-e8p8-28-1} & \input{tableau-e8p8-28-2} & \input{tableau-e8p8-28-3} & \input{tableau-e8p8-28-4}
\end{array}
$$

$$
\begin{array}{ccc}
\input{tableau-e8p8-28-5} & \input{tableau-e8p8-28-6} & \input{tableau-e8p8-28-7}
\end{array}
$$

\newpage

\section{$F_4$}

\subsection{$F_4/P_1$}


\noindent
Schubert cells in degree $0$:

\noindent
$
\begin{array}{rcl}
\sigma_{(0,0,0,0)} & = & 1\ \\ 
\end{array}
\vspace{.5cm}
$

\noindent
Schubert cells in degree $1$:

\noindent
$
\begin{array}{rcl}
\sigma_{(1,0,0,0)} & = & h\\ 
\end{array}
\vspace{.5cm}
$

\noindent
Schubert cells in degree $2$:

\noindent
$
\begin{array}{rcl}
\sigma_{(1,1,0,0)} & = & h^{2}\\ 
\end{array}
\vspace{.5cm}
$

\noindent
Schubert cells in degree $3$:

\noindent
$
\begin{array}{rcl}
\sigma_{(1,1,2,0)} & = & 1/2\ h^{3}\\ 
\end{array}
\vspace{.5cm}
$

\noindent
Schubert cells in degree $4$:

\noindent
$
\begin{array}{rcl}
\sigma_{(1,2,2,0)} & = & s\\ 
\end{array}
\\
\begin{array}{rcl}
\sigma_{(1,1,2,2)} & = & 1/4\ h^{4} - 1/2\ s\\ 
\end{array}
\vspace{.5cm}
$

\noindent
Schubert cells in degree $5$:

\noindent
$
\begin{array}{rcl}
\sigma_{(1,2,2,2)} & = & 1/4\ h^{5} - 1/2\ hs\\ 
\end{array}
\\
\begin{array}{rcl}
\sigma_{(2,2,2,0)} & = & -1/2\ h^{5} + 2\ hs\\ 
\end{array}
\vspace{.5cm}
$

\noindent
Schubert cells in degree $6$:

\noindent
$
\begin{array}{rcl}
\sigma_{(1,2,4,2)} & = & 1/4\ h^{6} - 3/4\ h^{2}s\\ 
\end{array}
\\
\begin{array}{rcl}
\sigma_{(2,2,2,2)} & = & -1/4\ h^{6} + h^{2}s\\ 
\end{array}
\vspace{.5cm}
$

\noindent
Schubert cells in degree $7$:

\noindent
$
\begin{array}{rcl}
\sigma_{(1,3,4,2)} & = & 3/8\ h^{7} - 5/4\ h^{3}s\\ 
\end{array}
\\
\begin{array}{rcl}
\sigma_{(2,2,4,2)} & = & -1/8\ h^{7} + 1/2\ h^{3}s\\ 
\end{array}
\vspace{.5cm}
$

\noindent
Schubert cells in degree $8$:

\noindent
$
\begin{array}{rcl}
\sigma_{(3,3,4,2)} & = & 1/4\ h^{8} - h^{4}s + 1/2\ s^{2}\\ 
\end{array}
\\
\begin{array}{rcl}
\sigma_{(2,4,4,2)} & = & -3/16\ h^{8} + 3/4\ h^{4}s - 1/4\ s^{2}\\ 
\end{array}
\vspace{.5cm}
$

\noindent
Schubert cells in degree $9$:

\noindent
$
\begin{array}{rcl}
\sigma_{(3,4,4,2)} & = & 3/16\ h^{9} - h^{5}s + 5/4\ hs^{2}\\ 
\end{array}
\\
\begin{array}{rcl}
\sigma_{(2,4,6,2)} & = & -3/16\ h^{9} + 7/8\ h^{5}s - 3/4\ hs^{2}\\ 
\end{array}
\vspace{.5cm}
$

\noindent
Schubert cells in degree $10$:

\noindent
$
\begin{array}{rcl}
\sigma_{(3,4,6,2)} & = & 1/16\ h^{10} - 1/2\ h^{6}s + h^{2}s^{2}\\ 
\end{array}
\\
\begin{array}{rcl}
\sigma_{(2,4,6,4)} & = & -1/8\ h^{10} + 11/16\ h^{6}s - 7/8\ h^{2}s^{2}\\ 
\end{array}
\vspace{.5cm}
$

\noindent
Schubert cells in degree $11$:

\noindent
$
\begin{array}{rcl}
\sigma_{(3,5,6,2)} & = & 9/32\ h^{11} - 15/8\ h^{7}s + 25/8\ h^{3}s^{2}\\ 
\end{array}
\\
\begin{array}{rcl}
\sigma_{(3,4,6,4)} & = & -1/8\ h^{11} + 11/16\ h^{7}s - 7/8\ h^{3}s^{2}\\ 
\end{array}
\vspace{.5cm}
$

\noindent
Schubert cells in degree $12$:

\noindent
$
\begin{array}{rcl}
\sigma_{(3,5,6,4)} & = & -1/8\ h^{8}s + 7/16\ h^{4}s^{2}\\ 
\end{array}
\vspace{.5cm}
$

\noindent
Schubert cells in degree $13$:

\noindent
$
\begin{array}{rcl}
\sigma_{(3,5,8,4)} & = & 1/16\ h^{13} - 15/32\ h^{9}s + 7/8\ h^{5}s^{2}\\ 
\end{array}
\vspace{.5cm}
$

\noindent
Schubert cells in degree $14$:

\noindent
$
\begin{array}{rcl}
\sigma_{(3,6,8,4)} & = & 33/128\ h^{14} - 27/16\ h^{10}s + 11/4\ h^{6}s^{2}\\ 
\end{array}
\vspace{.5cm}
$

\noindent
Schubert cells in degree $15$:

\noindent
$
\begin{array}{rcl}
\sigma_{(4,6,8,4)} & = & 75/128\ h^{15} - 119/32\ h^{11}s + 187/32\ h^{7}s^{2}\\ 
\end{array}
\vspace{.5cm}
$

\subsection{$F_4/P_4$}

\input{f4p4}

\newpage

\section{$G_2$}

\subsection{$G_2/P_1$}


\noindent
Schubert cells in degree $0$:

\noindent
$
\begin{array}{rcl}
\sigma_{(0,0)} & = & 1\ \\ 
\end{array}
\vspace{.5cm}
$

\noindent
Schubert cells in degree $1$:

\noindent
$
\begin{array}{rcl}
\sigma_{(1,0)} & = & h\\ 
\end{array}
\vspace{.5cm}
$

\noindent
Schubert cells in degree $2$:

\noindent
$
\begin{array}{rcl}
\sigma_{(1,1)} & = & h^{2}\\ 
\end{array}
\vspace{.5cm}
$

\noindent
Schubert cells in degree $3$:

\noindent
$
\begin{array}{rcl}
\sigma_{(3,1)} & = & 1/2\ h^{3}\\ 
\end{array}
\vspace{.5cm}
$

\noindent
Schubert cells in degree $4$:

\noindent
$
\begin{array}{rcl}
\sigma_{(3,2)} & = & 1/2\ h^{4}\\ 
\end{array}
\vspace{.5cm}
$

\noindent
Schubert cells in degree $5$:

\noindent
$
\begin{array}{rcl}
\sigma_{(4,2)} & = & 1/2\ h^{5} - q\\ 
\end{array}
\vspace{.5cm}
$

\subsection{$G_2/P_2$}


\noindent
Schubert cells in degree $0$:

\noindent
$
\begin{array}{rcl}
\sigma_{(0,0)} & = & 1\ \\ 
\end{array}
\vspace{.5cm}
$

\noindent
Schubert cells in degree $1$:

\noindent
$
\begin{array}{rcl}
\sigma_{(0,1)} & = & h\\ 
\end{array}
\vspace{.5cm}
$

\noindent
Schubert cells in degree $2$:

\noindent
$
\begin{array}{rcl}
\sigma_{(3,1)} & = & 1/3\ h^{2}\\ 
\end{array}
\vspace{.5cm}
$

\noindent
Schubert cells in degree $3$:

\noindent
$
\begin{array}{rcl}
\sigma_{(3,3)} & = & 1/6\ h^{3} - 1/2\ q\\ 
\end{array}
\vspace{.5cm}
$

\noindent
Schubert cells in degree $4$:

\noindent
$
\begin{array}{rcl}
\sigma_{(6,3)} & = & 1/18\ h^{4} - 1/2\ hq\\ 
\end{array}
\vspace{.5cm}
$

\noindent
Schubert cells in degree $5$:

\noindent
$
\begin{array}{rcl}
\sigma_{(6,4)} & = & 1/18\ h^{5} - 5/6\ h^{2}q\\ 
\end{array}
\vspace{.5cm}
$

\end{document}